\documentclass{amsart}
\usepackage{amssymb,latexsym,amsmath,amsfonts}
\usepackage{cmap}
\input xy
\xyoption{all}
\usepackage[active]{srcltx}

\newcommand{\cal}{\mathcal}

\makeatletter
\renewcommand{\subsection}{\@startsection{subsection}{2}{0mm}{-2mm}{-2mm}{\bf\normalsize}}

\def\sbsnt#1{\subsection{#1}}
\makeatother

\newtheorem{formula}{}[section]
\newtheorem{definition}[formula]{Definition}
\newtheorem{corollary}[formula]{Corollary}
\newtheorem{remark}[formula]{Remark}
\newtheorem{lemma}[formula]{Lemma}
\newtheorem{theorem}[formula]{Theorem}

\def\thrm{\begin{theorem}}
\def\thrml#1{\begin{theorem}\label{#1}}
\def\ethrm{\end{theorem}}
\def\rmrk{\begin{remark}}
\def\rmrkl#1{\begin{remark}\label{#1}}
\def\ermrk{\end{remark}}
\def\dfntn{\begin{definition}}
\def\dfntnl#1{\begin{definition}\label{#1}}
\def\edfntn{\end{definition}}
\def\nmrt{\begin{enumerate}}
\def\enmrt{\end{enumerate}}
\def\tm#1{\item[{\rm (#1)}]}
\def\qtn{\begin{equation}}
\def\qtnl#1{\begin{equation}\label{#1}}
\def\eqtn{\end{equation}}
\def\lmm{\begin{lemma}}
\def\lmml#1{\begin{lemma}\label{#1}}
\def\elmm{\end{lemma}}
\def\crllr{\begin{corollary}}
\def\crllrl#1{\begin{corollary}\label{#1}}
\def\ecrllr{\end{corollary}}
\def\css{\begin{cases}}
\def\ecss{\end{cases}}

\def\proof{\noindent{\bf Proof}.\ }

\def\P{{\cal P}}

\def\cS{{\cal S}}
\def\T{{\cal T}}
\def\W{{\cal W}}
\def\X{{\cal X}}
\def\Y{{\cal Y}}

\def\FF{{\mathbb F}}
\def\ZZ{{\mathbb Z}}

\def\SD{{\scriptscriptstyle\Delta}}
\def\SG{{\scriptscriptstyle\Gamma}}

\DeclareMathOperator{\aut}{Aut}

\DeclareMathOperator{\AGL}{AGL}
\DeclareMathOperator{\AGaL}{A{\rm \Gamma}L}

\DeclareMathOperator{\fis}{Fis}
\DeclareMathOperator{\Fit}{Fit}
\DeclareMathOperator{\Fix}{Fix}
\DeclareMathOperator{\fix}{fix}
\DeclareMathOperator{\GaL}{{\rm \Gamma}L}

\DeclareMathOperator{\GF}{GF}
\DeclareMathOperator{\GL}{GL}
\DeclareMathOperator{\id}{id}
\DeclareMathOperator{\img}{Im}

\DeclareMathOperator{\inv}{Inv}

\DeclareMathOperator{\iso}{Iso}

\DeclareMathOperator{\orb}{Orb}

\DeclareMathOperator{\spp}{Sp}
\DeclareMathOperator{\spa}{Span}

\DeclareMathOperator{\sym}{Sym}
\DeclareMathOperator{\Test}{Test}

\DeclareMathOperator{\Wr}{Wr}

\def\bul{\hfill\vrule height .9ex width .8ex depth -.1ex}
\def\bull{\hfill\vrule height .9ex width .8ex depth -.1ex\medskip}

\def\qaq{\quad\text{and}\quad}

\def\lg{\langle}

\def\pwr{\uparrow}
\def\rg{\rangle}
\newcommand{\rsc}[3]{#1\circ_{\scriptscriptstyle_{\hspace{-1pt}#2}}\hspace{-1pt}#3}

\def\sss{\scriptscriptstyle}

\begin{document}
\title[Bases of schurian antisymmetric coherent configurations]{Bases of schurian
antisymmetric coherent configurations and\\ isomorphism test for schurian tournaments}
\author{Ilya Ponomarenko}
\address{Steklov Institute of Mathematics at St. Petersburg, Russia}
\email{inp@pdmi.ras.ru}
\thanks{The work was partially supported by RFFI Grant 11-01-00760-a}
\date{}
\maketitle

\begin{abstract}
It is known that for any permutation group $G$ of odd order one can find a subset of 
the permuted set whose stabilizer in $G$ is trivial, and if $G$ is primitive, then also 
a base of size at most~$3$. Both of these results are generalized to the coherent
configuration of~$G$ (that is in this case a schurian antisymmetric coherent configuration).
This enables us to construct a polynomial-time algorithm for recognizing and isomorphism
testing of schurian tournaments (i.e. arc colored tournaments the coherent configurations 
of which are schurian).
\end{abstract}

\section{Introduction}
Let $\X$ be a coherent configuration (as for the background of coherent configurations 
we refer to Section~\ref{100511a} and~\cite{EP09}). A {\it base} of $\X$ is a point set
$\Delta$ such that the smallest fission of $\X$ in which all points of $\Delta$ 
are fibers, is the complete coherent configuration.\footnote{In survey \cite{CB} 
the name EP-base was used.} The minimal size of $\Delta$ is called the {\it base number} 
of $\X$ and is denoted by $b(\X)$. It is easily seen that $0\le b(\X)\le n-1$ where 
$n$ is the degree of~$\X$. Besides, given a permutation group $G$ denote by
$b(G)$ the smallest size of a base of $G$.\footnote{A base of a permutation group is a set of permuted points whose pointwise stabilizer is trivial.} Then
\qtnl{090511a}
b(G)\le b(\inv(G))
\eqtn
where $\inv(G)$ is the coherent configuration associated with $G$. A weaker upper bound 
for $b(G)$ enables us to estimate the maximal order of uniprimitive group as it was done in~\cite{B81}. It also follows from Theorem~0.2 of that paper that the base number of a nontrivial primitive coherent configuration of degree~$n$ is less than 
$4\sqrt{n}\log n$.\medskip

The equality in \eqref{090511a} is obviously attained when the group $G$ is trivial 
or symmetric. A nontrivial example of the equality was found in \cite{EP01ce} for
$G$ being the automorphism group of a cyclotomic scheme over finite field. On the other
hand, in general inequality~\eqref{090511a} is strict even for solvable groups:
if $G$ is a solvable $2$-transitive group of degree $n$, then $b(G)\le 4$ by 
\cite{S96}, but in this case $\inv(G)$ is trivial, and hence $b(\inv(G))=n-1$. 
In contrast to this example we prove here the following theorem.

\thrml{290311a}
Let $G$ be a primitive permutation group of odd order. Then the base number of
the coherent configuration $\inv(G)$ is at most~$3$.
\ethrm

As an immediate consequence of inequality \eqref{090511a} and Theorem~\ref{290311a} we
deduce that the base number of a primitive permutation group of odd order is at most~$3$
(this result have been earlier proved in \cite{Es91b}). In proving Theorem~\ref{290311a} 
we also get a generalization of the Gluck theorem that in any odd order permutation group
the stabilizer of some subset of the permuted set is trivial~\cite{G83}. Namely, 
a {\it generalized base} of a coherent configuration~$\X$ with the point set $\Omega$
is a set $\Pi\subset 2^\Omega$ such that the smallest fission of $\X$ in which any element
of~$\Pi$ is a union of some fibers of~$\X$, is the complete coherent configuration.
The minimal size of the set $\Pi$ is called the {\it generalized base number} of $\X$ and 
is denoted by $gb(\X)$. Again, it is easily seen that
\qtnl{050711a}
gb(G)\le gb(\inv(G))
\eqtn
where $gb(G)$ is the minimal size of the set $\Pi$ for which the intersection of
all $G_{\{\Delta\}}$ with $\Delta\in\Pi$ is trivial. 

\thrml{200411c}
Let $G$ be a permutation group of odd order. Then the generalized base number of
the coherent configuration $\inv(G)$ is at most~$1$.
\ethrm

A coherent configuration $\X$ is called {\it schurian} if there exists a permutation
group~$G$ such that $\X=\inv(G)$. Over all coherent configurations schurian ones are 
relatively rare in occurence. For example for infinitely many positive integers $n$ 
there are exponentially many antisymmetric coherent configuration of rank~$3$ and 
degree $n$; on the other hand, such a configuration is schurian if and only if it 
arises from the Payley tournament on $n$ vertices. In this paper we apply
Theorem~\ref{290311a} to get the following result.

\thrml{270411a}
Given an antisymmetric coherent configuration $\X$ on $n$ points one can test in 
time $n^{O(1)}$ whether $\X$ is schurian, and (if so) find the group $\aut(\X)$.
\ethrm

Antisymmetric coherent configurations are closely related with tournaments
(we recall that tournament is a directed graph in which any two distinct
vertices are joined by a unique arc). Indeed, one can easily seen that if
$(\Omega,S)$ is such a configuration and $A$ is a maximal subset of~$S$ such 
that $A\cap A^*=\emptyset$, then $(\Omega,A)$ is a tournament. Conversely,
the coherent configuration obtained from an arc colored tournament $T$ by means 
of the Weisfeiler-Leman algorithm\footnote{This algorithm was given in detail 
in \cite{W76}; see also Subsection~\ref{310711a}.} is antisymmetric. 
When this configuration is schurian, we say that the tournament 
$T$ is {\it schurian}. In particular, this is always the case when the color
classes of arcs are the orbits of the group $\aut(T)$ (acting on the pairs
of vertices).\medskip

Let us turn to the tournament isomorphism problem. It is a special case of the Graph
Isomorphism Problem that consists in finding an efficient algorithm to test whether 
or not two given (arc colored) tournaments are isomorphic. At present, the best result 
here is the algorithm from \cite{BL83} testing the isomorphism of two $n$-vertex 
tournaments in time $n^{O(\log n)}$ (see also \cite{ADM}). In this paper we prove
the following result.

\thrml{140211x}
Let $\T_n$ be the class of all schurian tournaments on $n$ vertices. Then the following
problems can be solved in time $n^{O(1)}$:
\nmrt
\tm{1} given a tournament $T$ on $n$ vertices, test whether $T\in\T_n$,
\tm{2} given a tournament $T\in\T_n$ find the group $\aut(T)$,
\tm{3} given tournaments $T_1,T_2\in\T_n$ find a set $\iso(T_1,T_2)$.
\enmrt
\ethrm

The proof of Theorem~\ref{140211x} is reduced to Theorem~\ref{270411a}). In the special 
case when $\X=\inv(G)$ for an odd order group $G$, the group $\aut(\X)$ is by definition
the $2$-closure of~$G$. Thus our algorithm, in particular, constructs in polynomial
time the $2$-closure of any odd order permutation group, that generalizes the main
result in~\cite{EP01}.\medskip

For the reader convenience we collect the basic facts on coherent configurations,
their bases and linear primitive solvable groups in Sections~\ref{100511a},
\ref{050711x} and~\ref{240411q} respectively. In Sections~\ref{050711q},
\ref{050711f} and~\ref{200311a} we prove some upper bounds for the numbers $b(\X)$ and
$gb(\X)$ when $\X$ is the wreath product or the exponentiation of coherent configurations. In Section~\ref{050711s} we give a sufficient condition for a
coherent configuration $\inv(G)$ with transitive $G$ to have a base of size at
most~$2$. This condition is used in Section~\ref{240511d}
where we prove that the equality in~\eqref{090511a} attained when $G$ is an affine
linear group with irreducible zero stabilizer (Theorem~\ref{291210d}). Finally, the
proofs of Theorems~\ref{290311a}, \ref{200411c} and~\ref{140211x} are given in
Section~\ref{050711u}.\medskip

{\bf Notation.}
Throughout the paper $\Omega$ denotes a finite set. The diagonal of the Cartesian
product $\Omega^2$ is denoted by~$1_\Omega$.\medskip

For $r\subset\Omega^2$ set
$r^*=\{(\beta,\alpha):\ (\alpha,\beta)\in r\}$. For $\Gamma,\Delta\subset\Omega$ set $r_{\SG,\SD}=r\cap(\Gamma\times\Delta)$ and $r_\SG=r_{\SG,\SG}$.\medskip

For any $\alpha\in\Omega$ set $1_\alpha=1_{\{\alpha\}}$
and $\alpha r=\{\beta\in\Omega:\ (\alpha,\beta)\in r\}$. \medskip

For $r,s\subset\Omega^2$ set
$r\cdot s=\{(\alpha,\beta)\in\Omega^2:\ (\alpha,\gamma)\in r,\ (\gamma,\beta)\in s$
for some $\gamma\in\Omega\}$, and set 
$r\otimes s=\{(\alpha,\beta)\in\Omega^2\times\Omega^2:\ (\alpha_1,\beta_1)\in r$
and $(\alpha_2,\beta_2)\in s\}$.\medskip 

For $S\in 2^{\Omega^2}$ denote by $S^\cup$
the set of all unions of the elements of $S$, and set $S^*=\{s^*:\ s\in S\}$
and $\alpha S=\cup_{s\in S}\alpha s$. For $T\in 2^{\Omega^2}$
set $S\cdot T=\{s\cdot t:\ s\in S,\, t\in T\}$.\medskip

For a permutation $g$ set $\fix(g)$ to be the number of points that $g$ leaves fixed. 
For a set $K$ of permutations set $\fix(K)=\max_g\fix(g)$ and $\Fix(K)=\sum_g\fix(g)$
where $g$ runs over the set $K^\#=K\setminus\{1\}$.

\section{Coherent configurations}\label{100511a}

Unfortunately up to now there is no commonly used terminology and notations in
the coherent configuration theory. In what follows we use a mix from~\cite{EP09}
and~\cite{Zi1}. All the facts presented below can be found in one of these sources.

\sbsnt{Definitions.}
A pair $\X=(\Omega,S)$ where $\Omega$ is a finite set and $S$ a partition of~$\Omega^2$, 
is called a {\it coherent configuration} on $\Omega$ if
$1_\Omega\in S^\cup$, $S^*=S$, and given $r,s,t\in S$, the number
$$
c_{rs}^t=|\alpha r\cap\beta s^*|
$$
does not depend on the choice of $(\alpha,\beta)\in t$. The elements of $\Omega$, $S$,
$S^\cup$ and the numbers $c_{rs}^t$ are called the {\it points}, the 
{\it basic relations}, the {\it relations} and the {\it intersection numbers} 
of~$\X$, respectively. For the intersection numbers the following equalities hold:
\qtnl{181110a}
c_{r^*s^*}^{t^*}=c_{sr}^t
\qaq 
|t|c_{rs}^{t^*}=|r|c_{st}^{r^*}=|s|c_{tr}^{s^*},\qquad r,s,t\in S.
\eqtn
The numbers $|\Omega|$ and $|S|$ are called the {\it degree} and the {\it rank} of $\X$. A unique basic
relation containing a pair $(\alpha,\beta)\in\Omega^2$ is denoted by $r_\X(\alpha,\beta)$ or $r(\alpha,\beta)$. The set of basic
relations contained in $r\cdot s$ with $r,s\in S^\cup$ is denoted by~$rs$.

\sbsnt{Fibers and homogeneity.}\label{180511a}
The point set $\Omega$ is a disjoint union of {\it fibers} which are the elements
of the set
$$
\Phi(\X)=\{\Gamma\subset\Omega:\ 1_\SG\in S\}
$$
Given a union $\Delta$ of fibers denote by $S_\SG$ the set of all nonempty relations 
$r_\SG$ with $r\in S$. Then $\X_\SG=(\Gamma,S_\SG)$ is a coherent configuration, called the 
{\it restriction} of~$\X$ to~$\Gamma$.\medskip

For any basic relation $r\in S$ there exist uniquely determined fibers $\Gamma$ and $\Delta$
such that $r\subset\Gamma\times\Delta$. The number $|\gamma r|=c_{ss^*}^t$ with
$t=1_\SG$, does not depend on $\gamma\in\Gamma$. It is called the {\it valency} 
of~$r$ and denoted~$n_r$. The maximum of all valences is denoted by $n_{max}$.\medskip

The coherent configuration $\X$ is called {\it homogeneous} or a {\it scheme} if
$1_\Omega\in S$. In this case $n_r=n_{r^*}$ and $|r|=nn_r$ for all $r\in S$ where
$n=|\Omega|$. Thus equalities in~\eqref{181110a} may be rewritten as follows:
\qtnl{150410a}
c_{r^*s^*}^{t^*}=c_{sr}^t\quad\text{and}\quad
n_tc_{rs}^{t^*}=n_rc_{st}^{r^*}=n_sc_{tr}^{s^*}.
\eqtn

\sbsnt{Equivalence relations.}
Let us define the {\it support} of a relation $r\subset\Omega^2$ to be the 
minimal set $\Gamma\subset\Omega$ such that $r\subset\Gamma^2$. 
Saying that $e\in S^\cup$ is an equivalence relation we mean that $e$ is an 
equivalence relation on its support; the set of classes of $e$ is denoted by $\Omega/e$.
According to~\cite[Subsection~3.2]{EPT} any such $e$ is a union of uniform equivalence
relations\footnote{An equivalence relation is uniform if all its classes have the
same size.} belonging to $S^\cup$ and having pairwise disjoint supports
This implies the following statement to be used in the proof of Corollary~\ref{190311a}.

\lmml{200311d}
Let $\X$ be a coherent configuration on $\Omega$, $e\in S^\cup$ is an equivalence relation
and $I\subset\Omega/e$. Suppose that no two classes of~$e$, one is in $I$ and another one
not in $I$, have the same size. Then the union of all elements of $I$ belongs to the set
$\Phi(\X)^\cup$.\bul
\elmm

Any coherent configuration has {\it trivial} equivalence relations in $S^\cup$: $1_\Omega$
and~$\Omega^2$. A homogeneous coherent configuration is called {\it primitive} if there 
are no other equivalence relations in $S^\cup$; otherwise it is called 
{\it imprimitive}.\medskip

Let $e\in S^\cup$ be an equivalence relation. Then given $\Gamma\in\Omega/e$ one can
construct the {\it restriction} of~$\X$ to~$\Gamma$ that is the coherent configuration 
$$
\X_\SG=(\Gamma,S_\SG)
$$ 
with $S_\Gamma$ as in Subsection~\ref{180511a}. The {\it quotient} of~$\X$ modulo~$e$ is
defined to be the coherent configuration 
$$
\X_{\Omega/e}=(\Omega/e,S_{\Omega/e})
$$ 
where $S_{\Omega/e}$ is the set of all nonempty relations of the form 
$\{(\Gamma,\Delta):\ s_{\SG,\SD}\ne\emptyset\}$ with $s\in S$.

\sbsnt{Fissions and fusions.}
There is a natural partial order\, $\le$\, on the set of all coherent configurations on
the set $\Omega$. Namely, given two coherent configurations $\X=(\Omega,S)$ and
$\X'=(\Omega,S')$ we set
$$
\X\le\X'\ \Leftrightarrow\ S^\cup\subset (S')^\cup.
$$
In this case $\X$ and $\X'$ are called respectively a {\it fusion} of $\X'$ and a
{\it fission} of $\X$. This order is preserved under taking the restriction to 
a set and the quotient modulo an equivalence. The minimal and maximal elements 
with respect to that order are the {\it trivial} and the {\it complete} coherent
configurations on~$\Omega$:
the basis relations of the former one are the reflexive relation $1_\Omega$ and 
(if $n>1$) its complement in $\Omega^2 $, whereas
the relations of the latter one are all binary relations on~$\Omega$.\medskip

Given two coherent configurations $\X_1$ and $\X_2$ on $\Omega$ there is
a uniquely determined coherent configuration $\X_1\cap\X_2$ also on $\Omega$,
the relation set of which is $(S_1)^\cup\cap(S_2)^\cup$ where $S_i$ is the set of basis
relations of~$\X_i$, $i=1,2$. This enables us to define the smallest fission
of a coherent configuration $\X$ on $\Omega$ containing a given set $\cS$ of
binary relations on $\Omega$ as follows:
$$
\fis(\X,\cS)=\bigcap_{\Y:\ \cS\subset T^\cup}\Y
$$
where $\Y=(\Omega,T)$ is a coherent configuration. In what follows we will omit $\X$
when it is the trivial coherent configuration. Besides, for $\Pi\subset 2^\Omega$
and $\Gamma\subset\Omega$ we define respectively the {\it $\Pi$-fission} and
{\it $\Gamma$-fission} of~$\X$ by
$$
\fis(\X,\Pi)=\fis(\X,\cS_\Pi)\quad\text{and}\quad\fis(\X,\Gamma)=\fis(\X,\Pi_\Gamma)
$$
where $\cS_\Pi=\{1_\SD:\ \Delta\in\Pi\}$ and $\Pi_\Gamma=\{\{\gamma\}:\ \gamma\in\Gamma\}$.
Sometimes we will also write $\X_{\alpha,\beta,\ldots}$ instead of
$\fis(\X,\{\alpha,\beta,\ldots\})$. One can see that any set in $\Pi^\cup$
is a union of fibers of the $\Pi$-fission of~$\X$. The following lemma
immediately follows from the definitions.

\lmml{260511a}
Let $\X=(\Omega,S)$ be a coherent configuration and $\alpha\in\Omega$. Then for
all $r,s,t\in S$ we have
$$
\alpha r\in(\Phi_\alpha)^\cup\quad\text{and}\quad
t_{\alpha r,\alpha s}\in (S_\alpha)^\cup
$$
where $\Phi_\alpha$ and $S_\alpha$ are the sets of fibers and basis relations of
the coherent configuration $\X_\alpha$. Moreover, $|\beta t_{\alpha r,\alpha s}|=c_{rt}^s$
for all $\beta\in\alpha r$.\bul
\elmm

\sbsnt{Isomorphisms and schurity.}\label{170511a}
Two coherent configurations $\X_1$ and $\X_2$ are called {\it isomorphic} if there exists
a bijection between their point sets that induces a bijection between their sets of basic
relations. Such a bijection is called an {\it isomorphism} between $\X_1$ and $\X_2$;
the set of all of them is denoted by $\iso(\X_1,\X_2)$.\medskip

The group of all isomorphisms of a coherent configuration $\X=(\Omega,S)$ to itself
contains a normal subgroup
$$
\aut(\X)=\{f\in\sym(\Omega):\ s^f=s,\ s\in S\}
$$
called the {\it automorphism group} of~$\X$ where 
$s^f=\{(\alpha^f,\beta^f):\ (\alpha,\beta)\in s\}$. 
Coversely, let $G$ be a permutation group on $\Omega$ and $S$ the set of orbits of 
the componentwise action of $G$ on~$\Omega^2$. Then $\inv(G):=(\Omega,S)$ is a 
coherent configuration; it is called the {\it coherent configuration of~$G$}. 
This coherent configuration is homogeneous if and only if the group $G$ is transitive.
One can also see that 
$$
\X\le\X'\ \Rightarrow\ \aut(X)\ge\aut(\X')\qaq
G\le G'\ \Rightarrow \inv(G)\ge\inv(G').
$$

A coherent configuration $\X$ is called {\it schurian} if $\X=\inv(G)$ for some 
permutation group $G$. In this case the group $G$ can be always replaced by~$\aut(\X)$.
Moreover, the schurity of $\X$ implies the schurity of all
its restrictions and quotients. An important example of a schurian scheme is a 
{\it cyclotomic scheme} over a finite field $\FF$; in this case
$G$ is an affine subgroup of $\AGL(1,\FF)$.~\footnote{In what follows saying that $G$ is an
affine (sub)group we mean that $G$ contains all the translations of the underlying linear
space.} In this paper we also deal with the scheme of a primitive solvable group.
The structure of such a group is given in the following statement proved
in~\cite[Section~4]{Su76}.

\thrml{140311c}
Let $G\le\sym(\Omega)$ be a primitive solvable permutation group. Then $|\Omega|=p^d$
for a prime $p$ and integer $d\ge 1$. Moreover the set $\Omega$ can be identified with a
linear space of dimension $d$ over field of order~$p$ so that
$$
G\le\AGL(d,p)\quad\text{and}\quad K\le\GL(d,p)
$$
where $K$ is the stabilizer of zero point in $G$; the group $G$ is 
affine and the group $K$ is irreducible.\bull
\ethrm

\sbsnt{Algebraic isomorphisms.}\label{130711a}
Let $\X=(\Omega,S)$ and $\X'=(\Omega',S')$ be coherent configurations. A bijection
$\varphi:S\to S',\ r\mapsto r'$ is called an {\it algebraic isomorphism} from~$\X$
to~$\X'$ if
\qtnl{f041103p1}
c_{r^{}s^{}}^{t^{}}=c_{r's'}^{t'},\qquad r,s,t\in S;
\eqtn
we say that $\X$ and $\X'$ are {\it algebraically isomorphic}. In this case they
have the same degree and rank. Moreover, $\varphi$ induces a bijection from $S^\cup$ onto
$(S')^\cup$ such that 
$$
(r\cup s)^\varphi=r^\varphi\cup s^\varphi,\qquad r,s\in S.
$$
This bijection preserves reflexive and equivalence relations. In particular,
we can define a bijection from $\Phi(\X)^\cup$ onto $\Phi(\X')^\cup$ so that 
$(1_\Gamma)^\varphi=1_{\Gamma^\varphi}$. Finally, given a set $\Gamma\in\Phi(\X)^\cup$
and an equivalence relation $e\in S^\cup$ we have the induced algebraic isomorphisms
$$
\varphi_\SG:\X^{}_{\SG^{}}\to\X'_{\SG'}\qaq
\varphi^{}_{\Omega^{}/e^{}}:\X^{}_{\Omega^{}/e^{}}\to \X'_{\Omega'/e'}
$$
where $\Gamma'=\Gamma^\varphi$ and $e'=e^\varphi$.\medskip

Any isomorphism~$f\in\iso(\X,\X')$ induces an algebraic isomorphism 
$r\mapsto r^f$ from~$\X$ to~$\X'$. The set of all isomorphisms 
inducing the algebraic isomorphism~$\varphi$ is denoted by $\iso(\X,\X',\varphi)$. 
Clearly,
$$
\iso(\X,\X,\id_S)=\aut(\X)
$$
where $\id_S$ is the identity on $S$. Let us give another example of algebraic isomorphism. 
Suppose that the scheme $\X$ is imprimitive  and $e\in S^\cup$ an equivalence relation. 
Then given any two sets $\Gamma,\Gamma'\in\Omega/e$ the mapping 
\qtnl{140711a}
\varphi_{\SG,\SG'}:\X_{\Gamma^{}}\to\X_{\Gamma'},\ s_{\SG^{}}\to s_{\SG'}
\eqtn
is an algebraic isomorphism (here $s$ runs over the set of all basis relations of
$\X$ that are contained in~$e$).

\sbsnt{Antisymmetric and $1$-regular coherent configurations.} A coherent configuration
$\X$ is called {\it antisymmetric} if 
$$
s\in S\qaq s=s^*\quad \Rightarrow\quad s\subset 1_\Omega,
$$
or equivalently if the cardinality of any basis relation of $\X$ is an odd number. The
latter condition implies that the valences of $\X$ are odd and that the 
coherent configuration $\inv(G)$ is antisymmetric if and only if $G$ is the group
of odd order. One can prove that the class of antisymmetric coherent configurations is
closed with respect to taking fissions, restrictions and quotients. In particular,
the automorphism group of antisymmetric coherent configuration has odd order.\medskip

A coherent configuration~$\X$ is called {\it $1$-regular} if it has a {\it regular}
point; by definition a point $\alpha\in\Omega$ is regular in~$\X$, if
\qtnl{040409b}
r\in S\ \Rightarrow\ |\alpha r|\le 1.
\eqtn
The set $\Gamma$ of all regular points is a union of fibers and any basic relation of the
coherent configuration $\X_\Gamma$ has valency~$1$. When $\Omega=\Gamma$, the coherent
configuration~$\X$ is called {\it semiregular}, and {\it regular} in homogeneous case.
Thus regular schemes are exactly {\it thin schemes} in the sense of~\cite{Zi1}. One can
also define such a scheme by the condition that any basis relation $r$ of it is {\it thin},
i.e. that 
$$
|\alpha r|\le 1\qaq|\alpha r^*|\le 1
$$
for all $\alpha\in\Omega$. We note that the set of all thin relations on the same set
is closed with respect to $*$ and $\cdot$.



\sbsnt{The Weisfeiler-Leman algorithm.}\label{310711a} 
From the algorithmic point of view a coherent configuration $\X$ on $n$ points is given by
the set $S$ of its basis relations. In this representation one can check in time $n^{O(1)}$
whether $\X$ is homogeneous or imprimitive. Moreover, within the same time one can list
the fibers of $\X$, and find a nontrivial equivalence relation $e\in S^\cup$ (if it exists)
as well as the quotient of~$\X$ modulo~$e$.\medskip

The well-known Weisfeiler-Leman algorithm is described in detail in~\cite[Section~B]{W76}. 
The input of it is a set $\cS$ of binary relations on a set
$\Omega$, and the output is the coherent configuration $\X=\fis(\cS)$. The running time of
the algorithm is polynomial in sizes of~$\cS$ and~$\Omega$. The canonical version of the
Weisfeiler-Leman algorithm have been studied in Section~M of the above book (under the
name  simultaneous stabilization), where in fact the following statement was proved.

\thrml{thbw}
Let $\cS_i$ be a set of $m$ binary relations on a set of size $n$, $i=1,2$.
Then given a bijection $\psi:\cS_1\to\cS_2$ one can check in time $mn^{O(1)}$ whether or
not there exists an algebraic isomorphism $\varphi:\fis(\cS_1)\to\fis(\cS_2)$ such that
$\varphi|_{\cS_1}=\psi$. Moreover, if $\varphi$ does exist, it can be found within the
same time.\bull
\ethrm

\section{Bases of a coherent configuration}\label{050711x}

\sbsnt{Generalized base.}\label{180811a}
A set $\Pi\subset 2^\Omega$ is called a {\it generalized base} of a coherent
configuration~$\X$ if the $\Pi$-fission of it is complete. When $\Pi$ consists of
singletons, we call it the {\it base} of~$\X$ and identify it with the 
corresponding subset of~$\Omega$. It is easily seen that $\Pi$ is a generalized base 
of any fission of~$\X$, and that any $\Pi'\subset 2^\Omega$ that contains $\Pi$ is a
generalized base of~$\X$. It is also clear that replacing some elements of~$\Pi$ by their
complements in~$\Omega$ produces a generalized base of~$\X$. The following simple statement
will be used in Section~\ref{200311a}.\medskip

\lmml{170311z}
Let $\X$ be a coherent configuration on $\Omega$, $\Pi$ a generalized base of~$\X$
and $\alpha\in\Omega$. Set $\X'=(\X_\alpha)_{\Omega'}$ and $\Pi'=\{\Gamma':\ \Gamma\in\Pi\}$
where $\Gamma'=\Gamma\setminus\{\alpha\}$ for all $\Gamma\subset\Omega$. Then
$\Pi'$ is a generalized base of~$\X'$.
\elmm
\proof Denote by $\Y$ the direct sum of the one-point coherent configuration on $\{\alpha\}$
and $\Pi'$-fission of $\X'$, i.e. the smallest coherent configuration on~$\Omega$ such that
$$
\{\alpha\}\in\Phi(\Y)\qaq\Y_{\Omega'}=\fis(\X',\Pi')
$$
Then obviously $\Y\ge\X_\alpha$ and $\Pi\subset\Phi(\Y)^\cup$. Therefore
$$
\Y\ge\fis(X_\alpha,\Pi)\ge\fis(\X,\Pi). 
$$
Since $\Pi$ is a generalized base of~$\X$, it follows that $\fis(\X,\Pi)$, and hence $\Y$,
is a complete coherent configuration. This implies that so is
$\Y_{\Omega'}=\fis(\X',\Pi')$. Thus $\Pi'$ is a generalized base of~$\X'$.\bul

\sbsnt{Generalized base number.} The smallest cardinality $gb(\X)$ (resp. by $b(\X)$) of 
a generalized base (resp. of a base) of the coherent configuration~$\X$ is called the 
{\it generalized base number} (resp. the {\it base number}) of~$\X$. Obviously,
\qtnl{240511b}
gb(\X)\le b(\X).
\eqtn
Since any fiber of $\X$ is a union of fibers in any its fission, we also have
\qtnl{200411e}
gb(\X)\le\max_{\Gamma\in\Phi}gb(\X_\SG)\quad\text{and}\quad
b(\X)\le\sum_{\Gamma\in\Phi}b(\X_\SG)
\eqtn
where $\Phi=\Phi(\X)$. Moreover, from the remark made in Subsection~\ref{180811a} it
immediately follows that
$$
\X'\ge\X\quad\Rightarrow\quad gb(\X')\le gb(\X)\qaq b(\X')\le b(\X).
$$  
It was observed in~\cite{EP01ce} that any regular point of
a coherent configuration forms a base of it. Thus $b(\X)\le 1$ for any $1$-regular coherent
configuration~$\X$. In the following statement we will use the fact that the equality
\qtnl{180811c}
b(\X)=b(\aut(\X))
\eqtn 
holds for any cyclotomic scheme~$\X$ (statement~(2) of \cite[Theorem~1.2]{EP01ce}).

\thrml{200411b}
Let $\X$ be an antisymmetric cyclotomic scheme over a finite field. Then $gb(\X)\le 1$
and $b(\X)\le 3$.
\ethrm
\proof We note that any antisymmetric scheme of degree $>1$ has rank at least~$3$. By
the hypothesis this implies that $\X$ is a proper cyclotomic scheme in the sense
of~\cite{EP01ce}. Therefore by the McConnel theorem (inclusion~(1) of this paper)
this implies that $\aut(\X)\le\AGaL(1,\FF)$ where $\FF$ is the underlying finite field. Thus 
due to \eqref{180811c} we have
$$
b(\X)=b(\aut(\X))\le b(\AGaL(1,\FF))\le 3.
$$ 
To prove that $gb(\X)\le 1$ set $b=b(\X)$. Without loss of generality we can assume that
$b=2$ or $b=3$. Denote by $\Y$ the $\Pi$-fission of $\X$ with $\Pi=\{B\}$
where $B=\{\alpha_0,\ldots,\alpha_{b-1}\}$ is a base of~$\X$. Then it suffices to 
verify that 
\qtnl{180811b}
B\not\in\Phi(\Y).
\eqtn 
Indeed, in this case the set $B$ must be the union of $b$ fibers 
which are singletons because the size of any fiber of antisymmetric configuration is of odd
cardinality. But then $\Y=\fis(\X,B)$ is the complete configuration. Thus $\Pi$ is a
generalized base of~$\X$ and we are done.\medskip

To prove~\eqref{180811b} suppose on the contrary that $B\in\Phi(\Y)$. Then $b=3$ because any
antisymmetric scheme, and hence $\Y_B$, has odd degree. However, up to isomorphism there
is a unique antisymmetric scheme of degree~$3$, namely, the scheme of a regular group
of order~$3$. This implies that $r_\Y(\alpha_0,\alpha_1)=r_\Y(\alpha_0,\alpha_2)^*$,
and hence
\qtnl{180811d}
r(\alpha_0,\alpha_1)=r(\alpha_0,\alpha_2)^*
\eqtn
On the other hand, by the transitivity of the group $\aut(\X)$ we can assume that
$\alpha_0=0_\FF$. Then it is easily seen that the set of fixed points of the 
two-point stabilizer $\aut(\X)_{\alpha_0,\alpha_1}$ is an additive subgroup of~$\FF$. 
So from~\eqref{180811c} it follows that the set $\{\alpha_0,\alpha_1,-\alpha_2\}$ is 
also a base of~$\X$. Thus without loss of generality we can assume that
$$
r(\alpha_0,\alpha_1)\ne r(\alpha_0,\alpha_2)^*.
$$
However, this contradicts~\eqref{180811d}.\bul

\sbsnt{Bases of size at most $2$.} A symmetric relation $s\in S^\cup$ is called
{\it connected} if any two distinct points in $\Omega$ are joined by a path in the graph
$(\Omega,s)$. It is well-known that a scheme $\X$ is primitive if and only
if any non-reflexive relation $s\cup s^*$, $s\in S$, is connected.

\thrml{140311a}
Let $\X$ be a coherent configuration and $s\in S^\cup$ a connected relation.
Suppose that for any point $\alpha\in\Omega$ the coherent configuration
$(\X_{\alpha})_{\alpha s}$ is semiregular. Then any pair of distinct points
in $s$ forms a base of~$\X$. In particular,
$b(\X)\le 2$.
\ethrm
\proof Without loss of generality we can assume that $s\cap 1_\Omega=\emptyset$. Let
$(\alpha,\beta)\in s$. Set $\Gamma=\{\gamma\in\Omega:\ \{\gamma\}\in\Phi(\X_{\alpha,\beta})\}$.
Then obviously $\alpha,\beta\in\Gamma$. Moreover, given $\gamma\in\Gamma$ we have
\qtnl{250511a}
\gamma s\subset\Gamma\quad\text{or}\quad \gamma s\cap\Gamma=\emptyset.
\eqtn
Indeed, suppose on the contrary that there exist points $\gamma\in\Gamma$ and 
$\gamma_1,\gamma_2\in\gamma s$ such that $\gamma_1\in\Gamma$ and $\gamma_2\not\in\Gamma$.
Since the coherent configuration $(\X_{\gamma})_{\gamma s}$ is semiregular this implies 
that $\{\gamma_2\}\in\Phi(\X_{\gamma,\gamma_1})$. However, the coherent configuration
$\X_{\gamma,\gamma_1}$ is a fusion of $\X_{\alpha,\beta}$
because $\gamma,\gamma_1\in\Gamma$. Therefore $\{\gamma_2\}\in\Phi(\X_{\alpha,\beta})$, 
and hence $\gamma_2\in\Gamma$. Contradiction.\medskip

Denote by $\Gamma_0$ the set of all points $\gamma\in\Gamma$ for which 
$\gamma s\subset\Gamma$. Then $\alpha\in\Gamma_0$ because $(\alpha,\beta)\in s$,
$\beta\in\Gamma$ and the coherent configuration $(\X_{\alpha})_{\alpha s}$ is semiregular.
By~\eqref{250511a} this implies that that $\Gamma_0$ is the connectivity component of 
the graph $(\Omega,s)$ that contains the vertex~$\alpha$. Since this graph
is connected, this implies that $\Gamma_0=\Omega$. Therefore $\Gamma=\Omega$. By the
definition of $\Gamma$ this means any fiber of the coherent configuration
$\X_{\alpha,\beta}$ is singleton, and hence this configuration is complete. Thus
$\{\alpha,\beta\}$ is a base of~$\X$.\bull

The following special statement will be used in the proof of Lemma~\ref{140411a}.
Below given a nonnegative integer $m$  and relations $r,s\in S$ we denote 
by $\rsc{r}{m}{s}$ the set of all $t\in r^*s$ such that $c_{rt}^s\le m$.

\lmml{160411a}
Let $\X$ be an antisymmetric primitive schurian scheme and $r\in S$ a
non-reflexive relation such that $|\rsc{r}{2}{r}\,\cup\,\rsc{r}{2}{r^*}|>2n_r/3$ and 
$\rsc{r}{2}{r^*}\ne\emptyset$. Then $b(\X)\le 2$.
\elmm
\proof By the hypothesis $s=r\cup r^*$ is a connected non-reflexive relation of~$\X$.
So by Theorem~\ref{140311a} it suffices to verify that the coherent 
configuration $\X_0=(\X_\alpha)_{\alpha s}$ is semiregular for all $\alpha\in\Omega$. 
However, from the schurity of~$\X$ it follows that
$$
\alpha r,\alpha r^*\in\Phi(\X_0).
$$
Denote by $S_1$ and $S_2$ the set of thin relations in $(S_0)_{\alpha r,\alpha r^*}$
and in $(S_0)_{\alpha r}$ respectively. Then one can see that the coherent
configuration~$\X_0$ is semiregular if and only if the following inequalities hold:
\qtnl{020811a}
|S_1|>0\qaq|S_2|>n_r/3.
\eqtn
Let $u\in \rsc{r}{2}{r^*}$. Since $n_{u^{}}=n_{u^*}$, from~\eqref{150410a} we obtain that
$u^*\in \rsc{r^*}{2}{r}$. This implies that any point in $\alpha r$ has at most two
neighbors in the relation $u'=u^{}_{\alpha r^*,\alpha r}$. On the other hand, the 
valences of $\X_0$ are odd. Thus by Lemma~\ref{260511a} the relation~$u'$ contains 
a relation from $S_1$. This proves the first inequality in~\eqref{020811a}. Let now 
$v\in\rsc{r}{2}{r}$. Then again any point in $\alpha r$ has at most two neighbors in 
$v'=v_{\alpha r^{},\alpha r^{}}$, and the above argument shows that $v$ contains a 
relation from $S_2$. Thus by the lemma hypothesis we have
$$
|S_1|+|S_2|\ge|\rsc{r}{2}{r}\,\cup\,\rsc{r}{2}{r^*}|>2n_r/3.
$$ 
Therefore either $S_1$ or $S_2$ contains more that $n_r/3$ elements. In the latter case 
the second inequality in~\eqref{020811a} is clear, whereas in the former case it
follows because $S_2$ contains a set $t\cdot S_1^*$ where $t$ is an arbitrary
element from~$S_1$.\bul

\sbsnt{Bases and isomorphisms.} 
The following statement shows how bases are used to find isomorphisms between coherent configurations.

\thrml{270411s}
Let $\X$ and $\X'$ be coherent configuration on $n$ points. Then given an algebraic 
isomorphism $\varphi:\X\to\X'$ all the elements in the set $\iso(\X,\X',\varphi)$ 
can be listed in time $(bn)^{O(b)}$ where $b=b(\X)$.
\ethrm
\proof By exhaustive search in time $n^{O(b)}$ one can find a size~$b$ base $B$ of $\X$.
Obviously, any isomorphism from $\X$ onto $\X'$ takes  it to a base of~$\X'$. Therefore
$$
\iso(\X,\X',\varphi)=\bigcup_{B'}\bigcup_g\iso_g(\X,\X',\varphi)
$$ 
where $B'$ runs over all size $b$ point sets of $\X'$, $g$ runs over all bijections
from~$B$ onto~$B'$, and $\iso_g(\X,\X',\varphi)$ consists of all $f\in\iso(\X,\X',\varphi)$
such that $f|_B=g$. Since there are at most $(bn)^{O(b)}$ possibilities for a pair $(B',g)$, 
only we need is to find in time $n^{O(1)}$ the set $\iso_g(\X,\X',\varphi)$ for fixed 
such a pair. To do this set 
$$
\cS=S\cup \{1_{\{\alpha\}}:\ \alpha\in B\}\qaq
\cS'=S'\cup \{1_{\{\alpha'\}}:\ \alpha'\in B'\}
$$ 
where $S$ and $S'$ are the
sets of basis relations of~$\X$ and~$\X'$ respectively. Then obviously the coherent
configurations
$$
\fis(\cS)=\fis(\X,B)\qaq\fis(\cS')=\fis(\X',B')
$$ 
are complete. So any algebraic isomorphism between them is induced by exactly one bijection between their fiber sets. Thus the required 
statement immediate follows from Theorem~\ref{thbw} for the bijection $\psi:\cS\to\cS'$ 
defined by the following conditions: $\psi|_S=\varphi$ and $\psi(1_{\{\alpha\}})=1_{\{\alpha^g\}}$ for all $\alpha\in B$.\bull

The following technical notion was introduced in \cite[Section~3.2]{EP03be}. Let $\X$
be a scheme and $e_0,e_1\in S^\cup$ two equivalence relations such that $e_0\subset e_1$. 
Set
$$
\Omega_0=\{\Gamma_1/e_0:\ \Gamma_1\in\Omega/e_1\}.
$$
By a {\it majorant} of the group $G=\aut(\X)$ with respect to the pair $(e_0,e_1)$ we mean a
permutation group $H$ on a set $\Delta$ together with a family of bijections
$f_\SG:\Gamma\to\Delta$ where $\Gamma\in\Omega_0$, such that
\qtnl{190811a}
(G^\Gamma)^{f_\SG}\le H
\eqtn
where $G^\Gamma=G^{\Gamma_1/e_0}$ is the permutation group induced by the natural action of
the setwise stabilizer~$G_{\{\Gamma_1\}}$ on the set $\Gamma$. 

\crllrl{290411a}
In the above notation let $\Gamma$ be an element of $\Omega_0$. Suppose that
\qtnl{290411b}
\iso(\X_{\Gamma^{}},\X_{\Gamma'},\varphi_{\SG,\SG'})\ne\emptyset
\quad\text{for all}\quad \Gamma'\in \Omega_0
\eqtn
where $\varphi_{\SG,\SG'}$ is the algebraic isomorphism~\eqref{140711a} for
$\X=\X_{\Omega/e_0}$. Then the group $\aut(\X_\Gamma)$ together with any family of bijections
$f_{\SG'}\in\iso(\X_{\Gamma^{}},\X_{\Gamma'},\varphi_{\SG,\SG'})$, 
$\Gamma'\in\Omega_0$, is a {\it majorant} of $\aut(\X)$ with respect to the
pair $(e_0,e_1)$. Moreover, it can be constructed in time $(bn)^{O(b)}$ where 
$n=|\Gamma|$ and $b=b(\X_\Gamma)$.
\ecrllr
\proof For any $\Gamma'\in\Omega_0$ we obviously have inclusion 
$G^{\Gamma'}\le \aut(\X_{\SG'})$. On the other hand, given a bijection
$f_{\SG'}\in\iso(\X_{\Gamma^{}},\X_{\Gamma'},\varphi_{\SG,\SG'})$, we have
$$
\aut(\X_{\SG'})^{f_{\SG'}}=\aut(\X_{\SG'}^{f_{\SG'}})=\aut(\X_{\SG}).
$$
Thus inclusion \eqref{190811a} holds for $\Delta=\Gamma$ and $H=\aut(\X_\Gamma)$. 
This proves the first statement of the lemma. The second one follows from
Theorem~\ref{270411s}.\bull

\section{Primitive linear groups of odd order}\label{240411q}

The structure of a finite solvable linear primitive group was studied in~\cite{Su76,MW}.
The following theorem is just a specialization of~\cite[Theorem~2.2]{YY10} for the
groups of odd order.

\thrml{020311a}
Let $K\le\GL(d,p)$ be a primitive group of odd order. Then every normal abelian
subgroup of~$K$ is cyclic and $K$ has a series $1<U\le F\le A\le K$ of
normal subgroups such that the following statements hold:
\nmrt
\tm{P1} $\spa(U)=\GF(p^a)$ where $a$ is a divisor of~$d$,
\tm{P2} $C_K(F)\le F\le\Fit(K)$ and $|F:U|=e^2$ for some integer~$e$
such that each prime divisor of $e$ divides $p^a-1$,
\tm{P3} $A=C_K(U)$ and $A/F$ is isomorphic to a completely reducible subgroup
of the group $\prod_{i=1}^m\spp(2n_i,p_i)$ where $p_i$ and $n_i$ are defined from the
prime power decomposition $e=\prod_{i=1}^mp_i^{n_i}$,
\tm{P4} $|K:A|$ divides $a$ and $ae$ divides~$d$.\bull
\enmrt
\ethrm

From statements (P1) and (P4) it follows that $|U|\le u_{a,p}$ and $|K:A|\le a_0$
where $u_{a,p}$ and $a_0$ are the maximal odd divisors of $p^a-1$ and $a$ respectively.
Thus
\qtnl{040211a}
|K|\le u_{a,p}\cdot e^2\cdot s_e\cdot a_0
\eqtn
where $s_e$ is the maximal order of the group $A/F$ for a fixed $e$ (see statement~(P3)).
The following two lemmas
collect some special facts on the group $K$ from Theorem~\ref{020311a} that
are contained in papers~\cite{Es91a,Es91b} or obtained by means of computer package
GAP~\cite{GAP}.

\lmml{150211j}
Let $e$ be the number from Theorem~\ref{020311a}. Then one of the
following statements hold:
\nmrt
\tm{1} $e=1$ and $K\le\GaL(1,p^d)$,
\tm{2} $e\in\{5,9,11,13\}$ and $s_5\le |K:F|=3$, $s_9,s_{11}\le 5$, $s_{13}\le 7$,
\tm{3} $e\ge 15$ is an odd integer and $s_e\le e^2/2$.
\enmrt
\elmm
\proof Suppose first that $e=1$. Then from~(P2) it follows that $U=F$ is a normal abelian
self-centralizing subgroup of~$K$. By \cite[Lemma~2.2]{MW} this implies that $F$ is
irreducible. Thus by~\cite[Theorem~2.1]{MW} we conclude that $K\le\GaL(1,p^d)$
which proves the second part of statement~(1). For $e\ge 15$ the required inequality
immediately follows from the fact that any completely reducible odd order subgroup of the
group $\spp(2n_i,p_i)$ has a regular orbit on the underground linear space (see
\cite[Theorem~A]{Es91a}). To deal with the case $1<e<15$ we start with some
observation.\medskip

Suppose that $e$ is an odd prime. We claim that the group $K/F$ has an irreducible
representation in $\GL(2,e)$. Indeed, by statement~(1) of \cite[Theorem~2.2]{YY10} the group 
$F$ is a central product of $U$ and a characteristic subgroup $E$ of $K$ that contains 
an extraspecial subgroup $E_0$ of order~$e^3$ and exponent~$e$. In particular, 
$Z=E\cap U$ is a central subgroup of~$F$, $E/Z\cong F/U$ and $|E_0\cap Z|\le e$. Therefore
by (P2) we have $|E:Z|=|F:U|=e^2$, and
\qtnl{200811a}
F/U\cong E_0/(E_0\cap Z)\cong\ZZ_e\times\ZZ_e.
\eqtn
However, by statement~(2) of \cite[Theorem~2.2]{YY10} the group $F/U$ is a completely
reducible $K/F$-module. Therefore $K/F$ has a representation in $\GL(2,e)$. Suppose that
this representation is not irreducible. Then one can find a group $E'>Z$
such that $|E:E'|=e$ and the group $E'/Z$ is $K/F$-invariant. But then
obviously $E'$ is a normal abelian subgroup of~$K$. Moreover from~\eqref{200811a}
it follows that $|E'\cap E_0|=e^2$. Therefore $E'$ contains an elementary abelian
subgroup of order~$e^2$. Thus $E'$ is normal abelian non-cyclic subgroup of~$K$,
which is impossible by Theorem~\ref{020311a}. The claim is proved.\medskip

Let now $1<e<15$. By means of GAP we find that (a) there are no
odd order irreducible subgroups in $\GL(2,e)$ for $e=3,7$, (b) the maximal order
of an irreducible odd order subgroup in $\GL(2,e)$ for $e=5,11,13$ equals respectively
to $3,15$ and~$21$, and (c) the irreducible subgroups in $\GL(2,11)$ of order $15$
and in $\GL(2,13)$ of order $21$ are not subgroups of $\spp(2,11)$ and $\spp(2,13)$
respectively. Thus the required statement immediately follows 
from the above claim unless $e=9$. In the remaining case the same argument as in the claim
shows that $K/F$ has a representation in $\GL(4,e)$. Since up to conjugacy the latter
group has a unique irreducible odd order subgroup and the order of it is $5$, it suffices
to verify that that representation is irreducible. Suppose that this is not true.
Then as in the above claim one can check that there is no $K/F$-invariant subgroup 
of $E/Z$ of order~$e$. Therefore such a subgroup has order~$e^2$. But this is 
impossible by statement~(a) with $e=3$.\bull

\lmml{200411z}
In the notation of Theorem~\ref{020311a} we have $\fix(K)\le p^{\lfloor 4d/9\rfloor}$.
Moreover, if $g$ is an element of~$K$ of prime order~$q$, then
\nmrt
\tm{1} $\fix(g)\le p^{\lfloor d/q\rfloor}$ for $g\in F$, 
\tm{2} $\fix(g)\le p^{\lfloor d/3\rfloor}$ for $g\not\in F$ and $q\ne 3$.
\enmrt
\elmm
\proof Follows from Lemma~1.3 of~\cite{Es91a} and the proof of it.\bull

\section{Bases of the wreath product}\label{050711q}
In this section we fix a coherent configuration $\X_i=(\Omega_i,S_i)$, $i=1,2$.
The {\it wreath product} $\X_1\wr\X_2$ can be defined as the smallest coherent
configuration $\X=(\Omega,S)$ with $\Omega=\Omega_1\times\Omega_2$ such that the 
set $S^\cup$ contains the equivalence relation~$e$ with classes
$\Omega_\alpha=\Omega_1\times\{\alpha\}$, $\alpha\in\Omega_2$, and
$$
(\X_{\Omega_\alpha})^{\pi_\alpha}=\X_1,\qquad \X^\pi=\X_2
$$
for all $\alpha$ where $\pi_\alpha:\Omega_\alpha\to\Omega_1$ and $\pi:\Omega\to\Omega_2$
are the natural projections. In particular, $\pi_\alpha\in\iso(\X_{\Omega_\alpha},\X_1)$ 
and $\X_{\Omega/e}=\X_2$. When the coherent configuration $\X$ is homogeneous, we have
\qtnl{270311a}
S=\{s_1\otimes 1_{\Omega_2}:\ s_1\in S_1\}\,\cup\,
\{\Omega_1^2\otimes s_2:\ s_2\in S_2,\ s_2\ne 1_{\Omega_2}\}.
\eqtn
Any imprimitive schurian scheme is isomorphic to a fission of the wreath product of
two smaller schemes. In general case the set $\Phi(\X)$ consists of all
sets $\Gamma_1\times\Gamma_2$ where $\Gamma_1\in\Phi(\X_1)$ and
$\Gamma_2\in\Phi(\X_2)$, and
\qtnl{240911a}
\X_{\Gamma_1\times\Gamma_2}=(\X_1)_{\Gamma_1}\wr(\X_1)_{\Gamma_2}.
\eqtn

\lmml{220711a}
Let $\X=\X_1\wr\X_2$ and $\Pi\subset 2^\Omega$. Suppose that 
\nmrt
\tm{1} $\Pi_\alpha=\{\Gamma\cap\Omega_\alpha:\ \Gamma\in\Pi\}$ is a generalized base of $\X_{\Omega_\alpha}$ for all $\alpha\in\Omega_2$,
\tm{2} $\Pi_{\Omega/e}=\{\Gamma^\pi:\ \Gamma\in\Pi\}$ is a generalized base of $\X_2$.
\enmrt
Then $\Pi$ is a generalized base of $\X$.
\elmm
\proof Set $\Y=\fis(\X,\Pi)$. Then obviously $\Y^\pi\ge\fis(\X_2,\Pi_2)$ is a complete
configuration by condition~(2). This implies that $\Omega_\alpha\in\Phi(\Y)^\cup$ for
all $\alpha\in\Omega_2$. It follows that $\Gamma\cap\Omega_\alpha\in\Phi(\Y)^\cup$
for all $\Gamma\in\Pi$. Therefore $\Y_{\Omega_\alpha}\ge \fis(\X_{\Omega_\alpha},\Pi_\alpha)$
is a complete configuration for all $\alpha$. Consequently, any fiber of $\Y$ is a singleton,
which means that the coherent configuration $\Y$ is complete. Thus 
$\Pi$ is a generalized base of $\X$.\bull

Let $\Pi$ be a generalized base of the coherent configuration $\X$. We say that 
$\Pi$ is {\it proper} if there exists a set $\Gamma\in\Pi$ such that
$\Gamma\cap\Omega_\alpha$ is a proper subset of $\Omega_\alpha$ for all $\alpha\in\Omega_2$.
Clearly, such a base can exist only if $|\Omega_2|>1$.  

\thrml{190311a}
Let $\X=\X_1\wr\X_2$ and $b=\max\{gb(\X_1),gb(\X_2)\}$. Suppose that $\X_1$
is antisymmetric. Then $gb(\X)\le b$. Moreover, if $b>0$, then there exists a proper
generalized base of~$\X$ of size~$b$. 
\ethrm
\proof Without loss of generality we can also assume that $|\Omega_1|>1$ and $b>0$, and
that the coherent configurations $\X_1$, $\X_2$, and hence $\X$, are homogeneous
(see the first inequality in~\eqref{200411e} and equality~\eqref{240911a}). 
Let $\Pi_i$ be a generalized base of $\X_i$ of size~$b$, $i=1,2$. The assumption
implies that the set $\Pi_1$ contains a proper subset of~$\Omega_1$. Let us choose a
bijection  $\Gamma_1\mapsto\Gamma_2$ from $\Pi_1$ onto $\Pi_2$, and denote by $\Pi$ 
the set of all
\qtnl{220711q}
\Gamma=\Gamma_1\times\Gamma_2\ \cup\ \Gamma'_1\times\Gamma'_2
\eqtn
with $\Gamma_1\in\Pi_1$ where $\Gamma'_i$ is the complement to $\Gamma_i$ in~$\Omega_i$,
$i=1,2$. Then $|\Pi|=b$. So it suffices to verify that $\Pi$ is a generalized base of~$\X$ 
(in this case $\Pi$ is proper because $\Pi_1$ contains a proper subset
of~$\Omega_1$).\medskip

One can see that conditions~(1) and~(2) of Lemma~\ref{220711a} are satisfied for 
the union of~$\Pi$ and $\Pi'=\{\Gamma_1\times\Gamma_2:\ \Gamma_1\in\Pi_1\}$. So 
by this lemma the union is a generalized base of~$\X$. Thus 
we have to verify only that $\fis(\X,\Pi\cup\Pi')$ is a fission of $\Y=\fis(\X,\Pi)$, 
or, equivalently, that
\qtnl{220711i}
\Gamma_1\times\Gamma_2\in\Phi(\Y)^\cup
\eqtn
for all $\Gamma_1\in\Pi_1$. To do this denote by $e'$ the equivalence relation on the
the set~$\Gamma$ such that $\Gamma/e'=I\cup I'$ with
$$
I=\{\Gamma\cap\Omega_\alpha:\ \alpha\in\Gamma_2\}\qaq 
I'=\{\Gamma\cap\Omega_\alpha:\ \alpha\in\Gamma'_2\}.
$$
Then from~\eqref{220711q} it follows that $e'=\Gamma^2\cap e$ is also a relation of~$\Y$. Besides, since $\X_1$
is antisymmetric, exactly one of the numbers $\Gamma_1$ and $\Gamma'_1$ is odd. This
implies that the hypothesis of Lemma~\ref{200311d} is satisfied for $\X=\Y$ and $e=e'$. By 
this lemma the union of all elements of $I$ belongs to the set $\Phi(\Y)^\cup$. Since
the union is obviously equal to $\Gamma_1\times\Gamma_2$, we conclude that \eqref{220711i}
holds.\bull

Let $\Pi$ be a proper generalized base of the coherent configuration $\X$. We say that 
$\Pi$ is {\it thin} if $|\Gamma\cap\Omega_\alpha|\le 1$ for all $\Gamma\in\Pi$ and
$\alpha\in\Omega_2$. The following statements will be used in Section~\ref{200311a} to 
estimate the base size of an exponentiation.

\thrml{190311ai}
Let $\X=\X_1\wr\X_2$. Suppose that $\X_1$ is antisymmetric.
Then~$\X$ has a thin generalized base of size 
$b=b_1+\max\{0,b_2-\lceil b_1/2\rceil\}$ where $b_1=b(\X_1)$ and $b_2=gb(\X_2)$.
\ethrm
\proof Let $\Pi_1$ be a base of~$\X_1$ of size $b_1$, and $\Pi_2$ a generalized base
of~$\X_2$ of size~$b_2$. Suppose first that $2b_2\ge b_1$. Then 
$b=\lfloor b_1/2\rfloor+b_2$. Without loss of generality we can assume that $b_1$ is 
even (otherwise we add an extra point to~$\Pi_1$). Let us fix
\begin{itemize}
\item a point $\delta\in\Omega_1$,
\item a decomposition $\Pi_1=B\cup B'$ into two disjoint sets of equal size,
\item a fixed point free involution $\beta\mapsto\beta'$ on $\Omega_1$ taking $B$ to $B'$,
\item an injection $B\to\Pi_2,\ \beta\mapsto\Gamma_\beta$; set $\Pi'_2$ to be the complement
to its image.
\end{itemize}
Denote by $\Pi$ the family of sets $\Gamma$ and $\Gamma'$ defined below for 
all $\beta\in B$, and sets $\{\delta\}\times\Gamma_2$ for all $\Gamma_2\in\Pi_2'$,
\qtnl{230711a}
\Gamma=\{\beta\}\times\Gamma_\beta\ \,\cup\ \,\{\beta'\}\times\Gamma'_\beta\qaq
\Gamma'=\{\beta'\}\times\Gamma_\beta\ \,\cup\ \,\{\beta\}\times\Gamma'_\beta.
\eqtn
Since $|B|=b_1/2$ and $|\Pi_2'|=b_2-b_1/2$, the family $\Pi$ is of size $b=b_1/2+b_2$.
To complete the proof we will verify that $\Pi$ is a generalized base
of~$\X$ (in this case $\Pi$ is thin just by the definition).\medskip

To prove that the coherent configuration $\Y=\fis(\X,\Pi)$ is complete, we note that
$\Phi(\Y)^\cup$ contains the sets $\Gamma^*=\Gamma\cup\Gamma'$ where $\Gamma$ and $\Gamma'$
are defined by~\eqref{230711a}. We claim that 
\qtnl{230711b}
\{\beta\}\times\Gamma_\beta\in\Phi(\Y)^\cup
\eqtn
for all $\beta\in B$. Then obviously $\{\beta'\}\times\Gamma'_\beta\in\Phi(\Y)^\cup$.
This implies that conditions~(1) and~(2) of Lemma~\ref{220711a} are satisfied for 
$\X$ and $\Pi^*$ where the latter consists of all sets
$\{\beta\}\times\Gamma_\beta$, $\{\beta'\}\times\Gamma'_\beta$ and
$\{\delta\}\times\Gamma_2$. So by this lemma $\Pi^*$ is a generalized base of~$\X$. Thus 
the coherent configuration $\Y\ge\fis(\X,\Pi^*)$  is complete and we are done.\medskip

To prove \eqref{230711b} suppose on the contrary that there is a set $\Delta\in\Phi(Y)$ 
such that
\qtnl{230711c}
\beta\alpha,\ \beta'\alpha'\in\Delta
\eqtn
for some $\beta\in B$, $\alpha\in\Gamma_\beta$ and $\alpha'\in\Gamma'_\beta$, 
where $\beta\alpha=(\beta,\alpha)$ and $\beta'\alpha'=(\beta',\alpha')$.
Denote by $e^*$ the equivalence relation on $\Gamma^*$ with classes $\Gamma^*\cap\Omega_\gamma=\{\beta,\beta'\}\times\{\gamma\}$ where
$\gamma\in\Omega_2$. Then $e^*=e\cap(\Gamma^*)^2$ is a relation of~$\Y$. 
Therefore the set $\Delta'=\Delta e^*$ belongs to $\Phi(\Y)^\cup$. Since the relation
$u:=e^*\setminus 1_{\Gamma^*}$ is thin, this implies that $\Delta'$ is a fiber of~$Y$
and $u_{\SD,\SD'}$ is a basic relation of~$\Y$. Thus from~\eqref{230711c} we obtain that
$$
r_\Y(\beta\alpha,\beta'\alpha)=r_\Y(\beta'\alpha',\beta\alpha')=u_{\SD,\SD'}.
$$
However, in this case $r_{\X_1}(\beta,\beta')=r_{\X_1}(\beta',\beta)$ which is impossible
because the coherent configuration $\X_1$ is antisymmetric.\medskip

Let now $2b_2<b_1$. Then $b=b_1$. In this case take two disjoint sets $B,B'\subset\Pi_1$ 
of the same size~$b_2$, choose a bijection from $B\to\Pi_2$, $\beta\mapsto\Gamma_\beta$, 
and set $\Pi_2'$ to be the family of $b_1-2b_2$ sets $\{\beta\}\times\Omega_2$ where
$\beta$ runs over the set $\Pi_1\setminus(B\cup B')$. Then the rest of the proof  is
completely analogous to the previous case.\bull

\section{Exponentiation}\label{050711f}
Let $\Gamma$ be a finite set, $m$ a positive integer and $\Delta=\{1,\ldots,m\}$.
Given a set $T\subset 2^{\Gamma\times\Gamma}$ denote by $T^{\otimes_m}$ the set of
all relations $t_1\otimes\cdots\otimes t_m$ with $t_i\in T$ for all~$i$. For a coherent
configuration $\Y=(\Gamma,T)$ the pair
$$
\Y^{\otimes_m}=(\Gamma^m,T^{\otimes_m})
$$
is also a coherent configuration (the Cartesian $m$-power of~$\Y$). Any permutation
group $L\le\sym(\Delta)$ has the natural action on $\Omega=\Gamma^m$: a permutation 
$l\in L$ moves a point $\alpha=(\ldots,\alpha_i,\ldots)$ to the point
$\alpha^l=(\ldots,\alpha_j,\ldots)$ with $j^l=i$ (and hence a relation 
$t=\cdots\otimes t_i\otimes\cdots$ to the relation $t^l=\cdots\otimes t_j\cdots\otimes$).
Denote by $T\pwr L$ the
set of all relations $t^L=\cup_{l\in L}t^l$ with $t\in T^{\otimes_m}$. Then the pair
\qtnl{160311c}
\X=\Y\pwr L=(\Omega,T\pwr L)
\eqtn
is a coherent configuration~\cite{EP99} called the {\it exponentiation} of~$\Y$
by~$L$.\footnote{It is a special case of the general construction of the exponentiation
introduced in~\cite{EP99}.} It was also proved in that paper that $\X$ is schurian if and
only if so is $\Y$, and that $\X$ is primitive if and only if $L$ is transitive and $Y$ 
is primitive and non-regular. It is easily seen that $\X=\Y$ whenever $m=1$.\medskip

In this paper we will use the exponentiation construction for the scheme of
a primitive solvable permutation group. The structure of such a group is described
in Theorem~\ref{140311c}. Depending on whether the group $K$ from this theorem is primitive 
(as a linear group) or not we will say that the scheme~$\X=\inv(G)$ is {\it linearly
primitive} or {\it linearly imprimitive}. In particular, in both cases $\X$ is schurian
and the following statement holds.

\thrml{160311a}
The scheme $\X$ has a (possibly trivial) fusion
isomorphic to $\Y\pwr L$ where $\Y$ is a linearly primitive scheme and $L$ is a transitive
group. Moreover, if $\X$ is antisymmetric, then $\Y$ is antisymmetric and~$L$ has odd order.
\ethrm
\proof Without loss of generality we can assume that the group $K$ is imprimitive. Then the linear space $\Omega$ is a direct sum
of the subspaces belonging the set
$$
\Delta=\{\Gamma^k:\ k\in K\}
$$
where $\Gamma$ is a proper
subspace of~$\Omega$, and $K$ is isomorphic to a subgroup of the wreath product of
the group $K^U=(K_{\{U\}})^U\le\GL(U)$ and the transitive permutation group
$K^\Delta\le\sym(\Delta)$
induced by the action of~$K$ on~$\Delta$, ~\cite[Section~15.2]{Su76}. According to
\cite[Proposition~4.1]{EP01} this implies that $G$ can be identified with a
subgroup of the wreath product $G^U\pwr K^\Delta$ of permutation groups
$G^U$ and $K^\Delta$ in primitive action. On the other hand, by \cite[p.212]{JKM} we have
$$
\inv(G^U\pwr K^\Delta)=\inv(G^U)\pwr K^\Delta
$$
Thus the scheme $\X=\inv(G)$ has a fusion $\Y\pwr L$ where $\Y=\inv(G^U)$ and $L=K^\Delta$.
Moreover, if the scheme $\Y$ is linear imprimitive, then by the above it has a fusion
$Y'\pwr L'$ for some scheme $\Y'=\inv(G')$ where $G'$ is a primitive solvable permutation
group and $L'$ is a transitive group. So by \cite[Proposition~3.3]{EP01} the scheme
$\X$ has a fusion $(\Y'\pwr L')\pwr L=\Y'\pwr(L'\wr L)$ and the first statement follows.
To prove the second statement it suffices to note that if the scheme $\X$ is antisymmetric, then
the group $K$ has odd order.\bull

\section{Bases of the exponentiation.}\label{200311a}
The following theorem gives upper bounds for the maximal sizes of generalized and ordinary
bases of the exponentiation~\eqref{160311c} when the coherent configuration~$\Y$
is antisymmetric. The former bound is the best possible whereas the latter one definitely not.
Nevertheless, even this rather weak bound is sufficient for the purpose of the paper.

\thrml{140311d}
Let $\Y$ be an antisymmetric coherent configuration and let $L$ be a transitive permutation 
group of odd order. Then $gb(\Y\pwr L)\le\max\{gb(\Y),b\}$ where $b=gb(\inv(L))$. Moreover, 
if $\Y$ is not complete, then
$$
b(\Y\pwr L)\le b(\Y)+\max\{0,b-\lceil (b(\Y)-1)/2\rceil\}.
$$
\ethrm

The proof of Theorem~\ref{140311d} will be given in the end of this section.
Let us fix some notations. Let $\X=(\Omega,S)$ be the coherent
configuration defined by~\eqref{160311c}. For any $i\in\{0,\ldots,m\}$
set
$$
r_i=\{(\alpha,\beta)\in\Gamma^m\times\Gamma^m:\ d(\alpha,\beta)=i\}
$$
where $d(\alpha,\beta)$ is the number of all $j\in\Delta$ such that $\alpha_j\ne\beta_j$. 
One can see that $r_i$ is the union of the relations from $T^{\otimes_m}$ in which $i$
factors are equal to $1_\Gamma$ and the other $m-i$ are $\Gamma^2\setminus 1_\Gamma$. 
Therefore $r_i\in S^\cup$ for all~$i$ (which means that $\X$ is a fission of a Hamming scheme).  
In what follows we set  $r_{-1}=\emptyset$ and $r=r_1$.\medskip

Let us fix a point $\gamma_0\in\Gamma$, and set $\alpha=\alpha(\gamma_0)$ to be the point 
of~$\Omega$ with all coordinates equal to~$\gamma_0$. Then the neighborhood $\alpha r$ 
of $\alpha$ in $r$ is the disjoint union of the sets
\qtnl{260311a}
\Gamma_i=\{\beta\in\Omega:\ d(\alpha,\beta)=1\ \text{and}\ \beta_i\ne\gamma\},\qquad i\in\Delta.
\eqtn
They are the classes of an equivalence relation on $\alpha r$ that is denoted by~$e$. It is
easily seen that $e=1_{\alpha r}\cup r_{\alpha r}$. Therefore $e$ is a relation of the
coherent configuration $\X_0=(\X_\alpha)_{\alpha r}$. The following two lemmas are key
ingredients in our proof.

\lmml{170311r}
The mapping $\rho:\Omega\to 2^{\alpha r}$, $\beta\mapsto\beta r_{d-1}\cap\alpha r$ 
where $d=d(\alpha,\beta)$, is an injection and
$$
\img(\rho)=\{\Lambda\subset\alpha r:\ |\Lambda\cap\Gamma_i|\le 1\ \text{for all}\ \,i\in\Delta\}.
$$
In particular, the set $\alpha r$ is a base of the coherent configuration~$\X_\alpha$.
\elmm
\proof Given $\beta\in\Omega$ and $i\in\Delta$ such that $\beta_i\ne\gamma_0$ set
$\beta^{(i)}$ to be the unique point in $\Gamma_i$ the $i$th coordinate of which is equal
to~$\beta_i$. Then obviously 
$$
d(\beta,\beta^{(i)})=d(\beta,\alpha)-1.
$$
Therefore $\beta^{(i)}\in\rho(\beta)$. On the other hand, let $\delta\in\rho(\beta)$.
Then $d(\delta,\beta)=d-1$. So the points $\delta$ and $\beta$ 
have exactly $m-d+1$ equal coordinates. At least $m-d$ of them equal $\gamma_0$. But
$\beta$ has exactly $m-d$ such coordinates. Therefore there is $i\in\Delta$ such that
$\beta_i\ne\gamma_0$ and $\beta_i=\delta_i$. This means that $\delta=\beta^{(i)}$. Thus
\qtnl{260311t}
\rho(\beta)=\{\beta^{(i)}:\ i\in \Delta,\ \beta_i\ne\gamma\}
\eqtn
which proves the first statement. To prove the second one it suffices to note that
no two points $\beta$ and $\beta'$ with $\rho(\beta)\ne\rho(\beta')$ belong the same
fiber of the coherent configuration $\fis(\X_\alpha,\alpha r)$.\bull

Set $\Gamma_0=\Gamma\setminus\{\gamma_0\}$. Let us define the mapping
$f:\Gamma_0\times\Delta\to\alpha r$ taking a pair $(\gamma,i)$ to the unique point
$\beta\in\Gamma_i$ for which $\beta_i=\gamma$. Then obviously $f$ is a bijection
and the $f$-image of the set $\Gamma_0\times\{i\}$ coincides with $\Gamma_i$ for 
all $i\in\Delta$.

\lmml{160311e}
Set $\Y_0$ to be the restriction of $\Y_{\gamma_0}$ to $\Gamma_0$. Then
$\X_0^{f^{-1}}\ge\Y_0\wr\inv(L)$.
\elmm
\proof Denote by $T_0$ the set of all relations $t_{\SG_0}$ with $t\in T$ (we recall 
that $T$ is the set of basis relations of~$\Y$). Then it is easily seen that
$\Y_0=\fis(T_0)$. So by the definition of wreath product it suffices to verify that 
for all $t_0\in T_0$ and all orbits $u\in\orb(L,\Delta^2)$ we have 
\qtnl{200711a}
(t_0\otimes 1_\Delta)^f,\ (\Gamma_0^2\otimes u)^f\in S_0^\cup
\eqtn
where $S_0$ is the set of basis relations of $\X_0$. To do this let $t_0\in T_0$. Then 
$t_0=t_{\SG_0}$  for some $t\in T$. By the definition of the exponentiation and the 
transitivity of $L$ the set $S$ contains the relation
\qtnl{200711b}
(t\otimes 1_\Delta\otimes\cdots\otimes 1_\Delta)^L=
(t\otimes 1_\Delta\otimes\cdots\otimes 1_\Delta)\,\cup\,\cdots\,\cup\,
(1_\Delta\otimes\cdots\otimes 1_\Delta\otimes t).
\eqtn
Denote by $s$ the restriction of this relation to $\alpha r$. Then $s\in S_0^\cup$ 
by Lemma~\ref{260511a}. On the other hand, given $i\in\Delta$ denote by $s_i$ the summand 
in the right-hand side of~\eqref{200711b} with $t$ being at the $i$th position. Then a
straightforward computation shows that $(s_i)_{\alpha r}$ coincides with the $f$-image
of $t_0\otimes 1_{\{i\}}$. It follows that the relation 
$$
(t_0\otimes 1_\Delta)^f=\bigcup_{i\in\Delta}(t_0\otimes 1_{\{i\}})^f=\bigcup_{i\in\Delta}s_i=s
$$
belongs to $S_0^\cup$ which proves the first part of~\eqref{200711a}. To prove
the second part let $u\in\orb(L,\Delta^2)$. Then $S^\cup$ contains
the union of relations $u_{ij}=s_1\otimes\cdots\otimes s_m$ with $s_i=u$, $s_j=u^*$ and
$s_k=1_\Gamma$ for all $k\ne i,j$. Denote by $s$ the restriction of this relation to $\alpha r$.
Then $s\in S_0^\cup$ by Lemma~\ref{260511a}. On the other hand, a straightforward
computation shows that given $(i,j)\in u$ the set $s_{ij}=(u_{ij})_{\alpha r}$ 
coincides with the $f$-image of the relation 
$(\gamma_0 u\times\{i\})\times(\gamma_0 u^*\times\{j\})$.
It follows that $e\cdot s_{ij}\cdot e=\Gamma_i^2\cup\Gamma_j^2\cup\Gamma_i\times\Gamma_j$.
Thus the relation
$$
(\Gamma_0^2\otimes u)^f=
\bigcup_{(i,j)\in u}( (\Gamma_0\times\{i\})\times(\Gamma_0\times\{j\}))^f=
(\bigcup_{(i,j)\in u}e\cdot s_{ij}\cdot e)\setminus e=(e\cdot s\cdot e)\setminus e
$$
belongs to $S_0^\cup$, and we are done.\bull

{\bf Proof of Theorem~\ref{140311d}.} To prove the first statement without loss of
generality we can assume that $b>0$, for otherwise, $|\Delta|=1$ and 
$\Y\pwr L=\Y$. Besides, by Lemma~\ref{170311z} we have $gb(\Y_0)\le gb(\Y)$
where $\Y_0$ is the coherent configuration defined in Lemma~\ref{160311e} with
arbitrarily chosen point~$\gamma_0$. Thus by Theorem~\ref{190311a} the coherent
configuration $\Y_0\wr\inv(L)$ has a proper generalized base $\Pi_0$ of size 
$$
b_0\le \max\{gb(\Y_0),b\}\le\max\{gb(\Y),b\}.
$$ 
By Lemma~\ref{160311e} this implies that the coherent configuration
$\X_0=(\X_{\alpha})_{\alpha r}$ with
$\alpha=\alpha(\gamma_0)$, has a generalized base $\Pi$ of size~$b_0$ that contains 
an element~$\Lambda_0$ such that
\qtnl{240711a}
0<|\Lambda_0\cap\Gamma_i|<|\Gamma_i|\quad\text{for all}\ i\in\Delta,
\eqtn 
where the sets $\Gamma_i$ are defined in \eqref{260311a}.
By the second statement of Lemma~\ref{170311r} the set $\Pi$ is a generalized 
base of the coherent configuration~$\X_\alpha$. Set $\Phi$ be the fiber of 
$\fis(\X,\Pi)$ that contains~$\alpha$. Then it suffices to verify that $\Phi=\{\alpha\}$
(indeed, in this case $\fis(\X,\Pi)\ge\fis(\X_\alpha,\Pi)$ and we are done). To do this
suppose that $\beta\in\Phi$. Then since $\alpha\in\Phi$, $\Lambda_0\subset\alpha r$ and 
$\Lambda_0$ is the union of fibers of $\fis(\X,\Pi)$, we have
\qtnl{210711s}
\Lambda_0\subset\beta r.
\eqtn
Then obviously $d(\alpha,\beta)\le 2$. So without loss of generality we can assume that 
\qtnl{210711w}
\alpha_1=\gamma\ne\beta_1\qaq\alpha_3=\gamma=\beta_3
\eqtn
(since $L$ is a transitive group of odd order, we can assume that $m\ge 3$). However,
by~\eqref{240711a} there exists a point $\delta\in\Lambda_0\cap\Gamma_3$. Then 
by~\eqref{210711w} we have $d(\delta,\beta)\ge 2$. So $\delta\not\in\beta r$ which
contradicts~\eqref{210711s}.\medskip

In the proof of the second statement we keep the notations of the previous paragraph.
Since the coherent configuration $\Y$ is not complete, we can choose the point 
$\gamma_0\in\Gamma$ so that there is a base of $\Y$ of size $b(\Y)$ that 
contains~ $\gamma_0$. Then by Lemma~\ref{170311z} the coherent configuration 
$\Y_0$ has a base of size at most $b(\Y)-1$. So by Theorem~\ref{190311ai} the coherent
configuration $\Y_0\wr\inv(L)$ has a thin generalized base of size 
\qtnl{240711c}
b_0\le (b(\Y)-1)+\max\{0,b-\lceil (b(\Y)-1)/2\rceil\}.
\eqtn
By Lemma~\ref{160311e} this implies that the coherent configuration $\X_0$ has a generalized
base $\Pi$ of size~$b_0$ such that $|\Lambda\cap\Gamma_i|\le 1$ for all $\Lambda\in\Pi$.
By Lemma~\ref{170311r} any such $\Lambda$ is of the form $\rho(\beta)$ for uniquely
determined point $\beta=\beta(\Lambda)$ in $\Omega$. Set 
$$
B_0=\{\beta(\Lambda):\ \Lambda\in\Pi\}.
$$
Then $\fis(\X_\alpha,B_0)\ge\fis(\X_\alpha,\Pi)$ because $\rho(\beta)$ is a union of fibers
of the coherent configuration $\X_{\alpha,\beta}$ for all $\beta\in B_0$. Since  $\Pi$ is a
generalized base $\X_0$, the second statement of Lemma~\ref{170311r} this implies that 
$B_0$ is a base of the coherent configuration~$\X_\alpha$. Thus the set $B=B_0\cup\{\alpha\}$
is a base of $\X$. Moreover, $|B|=|\Pi|+1=b_0+1$ and the required statement follows
from~\eqref{240711c}.\bull

\section{Indistinguishing number and base number}\label{050711s}

Let $\X=(\Omega,S)$ be a scheme. For any two points
$\alpha,\beta\in\Omega$ denote by $\Omega_{\alpha,\beta}$ the set of all $\gamma\in\Omega$
such that $r(\alpha,\gamma)=r(\beta,\gamma)$. Then
\qtnl{120211d}
|\Omega_{\alpha,\beta}|=\sum_{t\in S}c_{tt^*}^s
\eqtn
where $s=r(\alpha,\beta)$. It follows that this number does not depend on the choice
of $(\alpha,\beta)\in s$ and is denoted by $c(s)$; in~\cite{MP09}
it was called the {\it indistinguishing number} of~$s$. The maximal indistinguishing number
of a non-reflexive basis relation of~$\X$ is denoted by $c=c(\X)$. It is easily seen
that $c(\X)\ge 0$ and the equality is attained if and only if the scheme~$\X$ is regular.\medskip

The number $n-c$ where $n=|\Omega|$, was called in paper~\cite{B81} the distinguishing
number of the coherent configuration~$\X$. It was proved there that if $\X$ is primitive
and $|S|\ge 3$, then $b(\X)\le 4\sqrt{n}\log{n}$. In the following theorem we are interested 
in the base number when $\X$ is not necessarily
primitive and $c$ is rather small.

\thrml{120211a}
Let $\X$ be a scheme such that $4c(m-1)<n$ where $m=n_{max}$. Then the coherent configuration
$\X_\alpha$ is $1$-regular for any $\alpha\in\Omega$. In particular, $b(\X)\le 2$.
\ethrm
\proof Without loss of generality we can assume that $m\ge 2$ (otherwise the scheme
$\X$ is regular, and the statement is obvious). Let $\alpha\in\Omega$
and $r\in S$. Given $\beta\in\Omega$ denote by $\Omega_\beta$ the set of all pairs
$(\delta,\gamma)\in\alpha r\times\alpha r$ such that $\delta\ne\gamma$ and
$\beta\in\Omega_{\delta,\gamma}$. Then it is easily seen that
$$
|\Omega_\beta|=\sum_{s\in S}c_{rs}^t(c_{rs}^t-1)=
\sum_{s\not\in r\circ t}c_{rs}^t=\sum_{s\in S}c_{rs}^t-|r\circ t|=n_r-|r\circ t|
$$
where $t=r(\alpha,\beta)$ and $r\circ t=\{s\in r^*t:\ c_{rs}^t=1\}$. 
This implies that if $\beta\in \alpha S'_r$ where $S'_r$ is the set of all $t'\in S$ with $|r\circ t'|< n_r/2$, then the set $\Omega_\beta$ has at least $n_r/2$ 
elements. Therefore
\qtnl{260711a}
|\alpha S'_r|\cdot\frac{n_r}{2}\le
\sum_{\beta\in\alpha S'_r}|\Omega_\beta|\le
|T|
\eqtn
where $T$ is the union of all $\Omega_\beta$ with $\beta\in\alpha S'_r$. However, for each
pair $(\delta,\gamma)\in\alpha r\times\alpha r$ with $\delta\ne\gamma$ there are at
most $c$ points $\beta$ such that $(\delta,\gamma)\in\Omega_\beta$. So the set $T$ has
at most $n_r(n_r-1)c$ elements. By \eqref{260711a} and the lemma
hypothesis this implies that $|\alpha S'_r|\le 2(n_r-1)c\le 2(m-1)c<n/2$. Thus
\qtnl{140211bcd}
|\alpha S_r|=n-|\alpha S'_r|>\frac{n}{2}
\eqtn
where $S_r=\{t\in S:\ |r\circ t|>n_r/2\}$ is the complement to $S'_r$.\medskip

To complete the proof we will show that any $\beta\in\Omega$ for which the relation
$r=r(\alpha,\beta)$ is of valency~$m$, is a regular point of the coherent
configuration~$\X_\alpha$, i.e. that
\qtnl{140211g}
\beta\,r_{\scriptscriptstyle \X_\alpha}(\beta,\gamma)=\{\gamma\}
\eqtn
for all $\gamma\in\Omega$. To do this set $u=r(\alpha,\gamma)$.
Then inequality~\eqref{140211bcd} implies that $|\alpha S_r|>n/2$ and $|\alpha S_u|>n/2$.
Therefore the sets $S_r$ and $S_u$ contain a common relation, say~$v$.
It follows that neither $r\circ v$ nor $u\circ v$ is empty; take $s_\beta\in r\circ v$ and
$s_\gamma\in u\circ v$. Then by the definition of $\circ$ one can find points $\beta'$ and
$\gamma'$ in $\alpha v$ such that
\qtnl{170211z}
\beta' s_\beta^*\cap\alpha r=\{\beta\}\quad\text{and}\quad
\gamma' s_\gamma^*\cap\alpha u=\{\gamma\}.
\eqtn
Moreover, we have $|r\circ v|>n_r/2$ because $v\in S_r$. Therefore one can find two relations
$t_\beta$ and $t_\gamma$ in $r\circ v$ such that
\qtnl{170211y}
\beta' t_\beta^*\cap\alpha r=\{\delta\}=\gamma' t_\gamma^*\cap\alpha r
\eqtn
for some point $\delta\in\alpha r$. The obtained configuration is represented at
Fig.~\ref{f8}.
\def\VRT#1{*=<5mm>[o][F-]{#1}}
\begin{figure}[h]
$\xymatrix@R=20pt@C=30pt@M=0pt@L=5pt{
& \VRT{\delta} \ar@{->}^{t_\gamma}[rr]& & \VRT{\gamma'}\ar@{<-}^{s_\gamma}[rr] & & \VRT{\gamma}          & \\
& & & \VRT{\alpha} \ar@{->}_{r}[ull] \ar@{->}^{r}[dll] \ar@{->}_{v}[u] \ar@{->}^{v}[d] \ar@{->}^{u}[urr] & & & \\
& \VRT{\beta}\ar@{->}_{s_\beta}[rr]    & &\VRT{\beta'}\ar@{<-}^{t_\beta}[uull] & &                       & \\
}$
\caption{}\label{f8}
\end{figure}
By Lemma~\ref{260511a} the set $(S_\alpha)^\cup$ contains the relations
$$
a_1=(s_\beta)_{\alpha r,\alpha v},\quad a_2=(t^*_\beta)_{\alpha v,\alpha r},\quad 
a_3=(t_\gamma)_{\alpha r,\alpha v},\quad a_4=(s^*_\gamma)_{\alpha v,\alpha u},
$$
and hence the relation $a=a_1\cdot a_2\cdot a_3\cdot a_4$. On the
other hand, since $n_r=m$ and $v\in S_r$, from~\eqref{150410a} it follows that $n_v=m$.
This implies that $s_\beta^*,t_\beta^*,t_\gamma^*\in S(v,r)$. Therefore due to~\eqref{170211z}
and \eqref{170211y} we obtain that
\qtnl{180211a}
\beta\,a_1=\{\beta'\},\quad \beta'a_2=\{\delta\},\quad
\delta\, a_3=\{\gamma'\},\quad \gamma' a_4=\{\gamma\}.
\eqtn
Thus $\beta\,r_{\scriptscriptstyle\X_\alpha}(\beta,\gamma)\subset \beta a=\{\gamma\}$ 
whence \eqref{140211g} follows.\bull

\crllrl{120211b}
Let $G\le\AGL(\Omega)$ be an affine group acting on a linear space~$\Omega$ over a finite
field, and $K$ a one point stabilizer of~$G$. Suppose that
\qtnl{270711a}
4(k-1)\Fix(K)<n.
\eqtn
where $n=|\Omega|$ and $k=|K|$. Then $b(\inv(G))\le 2$. In particular, this is always true
whenever $4k(k-1)f<n$ where $f=\fix(K)$.
\ecrllr
\proof Set $\X=\inv(G)$. Choose two points $\alpha$ and $\beta$ such that
$c=c(s)=|\Omega_{\alpha,\beta}|$ where $s=r(\alpha,\beta)$. Then any point in $\Omega_{\alpha,\beta}$ is a fixed point of a permutation from the set 
$G_{\alpha\to\beta}=\{a\in G:\ \alpha^a=\beta\}$. Therefore 
\qtnl{140211e}
c=|\Omega_{\alpha,\beta}|\le\Fix(G_{\alpha\to\beta}).
\eqtn
On the other hand, any $a\in G_{\alpha\to\beta}$ is an affine mapping on~$\Omega$, say
$x\mapsto hx+b$ where $h$ is the matrix from $G_\alpha$ and $b=\beta-\alpha$ is a vector.
Therefore the numbers $\fix(a)$ and $\fix(h)$ are equal respectively to the numbers of
solutions of linear equation systems $(h-e)x=b$ and $(h-e)x=0$ where $e$ is the identity
matrix. When $0<\fix(a)\le n$, the latter numbers are equal. Therefore the right-hand
side of~\eqref{140211e} coincides with $\Fix(G_\alpha)$. 
Thus since $G_\alpha$ and $K$ are conjugate in $G$, we obtain from~\eqref{140211e} that
\qtnl{270711b}
c\le\Fix(K).
\eqtn
Next, since $\X$ is the scheme of the transitive group~$G$, we have $n_s\le k$ for all $s\in S$. 
Therefore it follows from~\eqref{270711b} and~\eqref{270711a} that
$$
4(m-1)c\le 4(k-1)\Fix(K)<n.
$$
Thus the first statement of the theorem follows from Theorem~\ref{120211a}. The
second statement also follows from the above inequality because $\fix(K)\le (k-1)f$.\bull

\section{Base of linearly primitive antisymmetric scheme}\label{240511d}
In this section we will prove that, in fact, the base number of a linearly primitive antisymmetric
scheme coincides with the base number of its automorphism group. However, our proof do not
use the fact that the latter number is at most~$3$.

\thrml{291210d}
The base number of a linearly primitive antisymmetric scheme is at most~$3$ and the
equality is attained only for cyclotomic schemes over a finite field.
\ethrm

In what follows we fix an affine group $G\le\AGL(d,p)$ such that the scheme $\X=\inv(G)$
is antisymmetric and linearly primitive. Then by Theorem~\ref{140311c} the zero stabilizer 
in~$G$ is an irreducible primitive group $K\le\GL(d,p)$ of odd order. For this group we keep 
the notation of Theorem~\ref{020311a}. In the following three lemmas we will verify that
\qtnl{220411a}
e>1\ \Rightarrow\ b(\X)\le 2.
\eqtn
In each of these lemmas we will subsequently exclude the values of $e$ for which the
implication could be violate, by means of Lemma~\ref{160411a} and Corollary~\ref{120211b}.

\lmml{050311a}
The implication \eqref{220411a} holds unless the quadruple $(e,a,d,p)$ is one of the
following:
\nmrt
\tm{E1} $(e,a,d)=(9,1,9)$ and $p\in\{7,13,19,31,37,43\}$,
\tm{E2} $(e,a,d,p)=(5,4,20,3),\ (5,3,15,11)$ or $(5,2,10,11)$,
\tm{E3} $(e,a,d)=(5,1,5)$ and $p\in\{11,\ldots,5591\}$, $p=1\mod 5$.
\enmrt
\elmm
\proof Suppose that the parameters $e,a,d,p$ of the group $K$ do not form a quadruple
from the lemma statement. Then by Corollary~\ref{120211b} it suffices to verify that
$p^d>4k^2f$ where $k=|K|$ and $f=\fix(K)$. However, $f\le p^{\lfloor 4d/9\rfloor}$
by Theorem~\ref{200411z}. Therefore due to~\eqref{040211a} the required inequality is a
consequence of the following one:
\qtnl{060211b}
p^d>4\cdot (u_{a,p}\cdot e^2\cdot s_e\cdot a_0)^2\cdot p^{\lfloor 4d/9\rfloor}.
\eqtn
Here $ae\le d$ by statement~(P4) of Theorem~\ref{020311a}. Therefore $a_0\le a\le d/e$
and $2u_{a,p}\le p^a\le p^{d/e}$. Besides, by the second and the third statements of
Lemma~\ref{150211j} we have $s_e\le e^2/2$. Consequently,
$4\cdot u_{a,p}\cdot s_e\cdot a_0\le p^{d/e}\cdot d\cdot e$. Thus to check
inequality~\eqref{060211b} it suffices to verify that
\qtnl{080211a}
4\cdot p^{d-\lfloor 4d/9\rfloor}>p^{2d/e}\cdot d^2e^6.
\eqtn
A direct computation shows that $4\cdot 3^{14d/27}>d^8$ for all $d\ge 54$. Therefore for all
integers $e\ge 54$ and all primes $p\ge 3$ the inequality
$$
4\cdot p^{d-\lfloor 4d/9\rfloor-2d/e}\ge 4\cdot 3^{(5/9-2/54)d}=4\cdot 3^{14d/27}>d^8\ge d^2e^6.
$$
holds for all $d\ge e$. This proves the required statement for all $e\ge 54$.\medskip

Denote by $d(e,p)$ the minimal positive integer~$d$ for which inequality~\eqref{080211a}
holds for a fixed~$e$ and~$p$, and by $p(a,e)$ the minimal element in the set $P(a,e)$ of
all odd primes $q$ such that each prime divisor of~$e$ divides $q^a-1$. Then by
statements~(P2) and~(P4) of Theorem~\ref{020311a} and by Lemma~\ref{150211j} without 
loss of generality we can assume that
\nmrt
\tm{C1} $e\in\{5,\ldots,53\}$ is an odd integer other than~$7$,
\tm{C2} when $e$ is fixed, $a\in\{1,\ldots,\lfloor d_0/e\rfloor\}$ where $d_0=d(e,3)$,
\tm{C3} when $e$ and $a$ are fixed, $p\in P(a,e)$.
\enmrt
For each $e$ satisfying (C1) we list in the Table~\ref{tbl1} below
the values of the function $d_3=d(e,3)$ (the second row), the possible values for the
integer~$a$ (the third row), and for a fixed $a$ also the values of the functions $p_a=p(a,e)$
and $d_p=d(e,p_a)$ (the fourth and the fifth rows respectively).\footnote{
When $e=5$, we have $p_a=3$ and $d(e,p_a)=55$ for $a=0\mod 4$, and $p_a=11$ and $d(e,p_a)=43$ 
otherwise.}
\def\sss#1{{$\scriptstyle #1$}}
\def\mbf{\mathbf}
\begin{table}[tbh]
\caption{}\label{tbl1}
\begin{center}
\begin{tabular}{|c|c|c|c|c| c|c|c|c|c| c|c|c|c|c| c|}
\hline
\sss{e}        & \sss{53} & \sss{51}  & \sss{49} & \sss{47}  & \sss{45} & \sss{43}  & \sss{41} & \sss{39} & \sss{37}  & \sss{35} & \sss{33} & \sss{31}  & \sss{29} & \sss{27} & \sss{25}  \\
\hline
\sss{d_3}      & \sss{54} & \sss{54}  & \sss{53} & \sss{53}  & \sss{53} & \sss{53}  & \sss{52} & \sss{52} & \sss{52}  & \sss{51} & \sss{51} & \sss{51}  & \sss{50} & \sss{50} & \sss{50}  \\
\hline
\sss{a}        & \sss{1}  & \sss{1}   & \sss{1}  & \sss{1}   & \sss{1}  & \sss{1}   & \sss{1}  & \sss{1}  & \sss{1}   & \sss{1}  & \sss{1}  & \sss{1}   & \sss{1}  & \sss{1}  & \sss{1,2} \\
\hline
\sss{p_a}      & \sss{107}& \sss{103} & \sss{29} & \sss{283} & \sss{31} & \sss{173} & \sss{83} & \sss{79} & \sss{149} & \sss{71} & \sss{67} & \sss{311} & \sss{59} & \sss{7}  & \sss{11}  \\
\hline
\sss{d_p}      & \sss{12} & \sss{12}  & \sss{16} & \sss{10}  & \sss{16} & \sss{10}  & \sss{12} & \sss{12} & \sss{10}  & \sss{12} & \sss{12} & \sss{9}   & \sss{12} & \sss{27} & \sss{22}  \\
\hline
\hline
\sss{e}        & \sss{23} & \multicolumn{2}{c|}{\sss{21}} & \multicolumn{2}{c|}{\sss{19}} & \multicolumn{2}{c|}{\sss{17}} & \multicolumn{2}{c|}{\sss{15}} & \multicolumn{2}{c|}{\sss{13}} & \sss{11}             & \multicolumn{2}{c|}{\sss{9}}                    & \sss{5}\\
\hline
\sss{d_3}      & \sss{49} & \multicolumn{2}{c|}{\sss{49}} & \multicolumn{2}{c|}{\sss{49}} & \multicolumn{2}{c|}{\sss{49}} & \multicolumn{2}{c|}{\sss{49}} & \multicolumn{2}{c|}{\sss{50}} & \sss{51}             & \multicolumn{2}{c|}{\sss{55}}                   & \sss{103} \\
\hline
\sss{a}        & \sss{1,2}& \sss{1}  & \sss{2}            & \sss{1}   & \sss{2}           & \sss{1}   & \sss{2}           & \sss{1,3} & \sss{2}           & \sss{1,2} & \sss{\mbf 3}      & \sss{{\mbf 1},2,3,4} & \sss{{\mbf 1},{\mbf 3},5} &  \sss{{\mbf 2},4,6} & \sss{1..20}\\
\hline
\sss{p_a}      & \sss{47} & \sss{43} & \sss{13}           & \sss{191} & \sss{37}          & \sss{103} & \sss{67}          & \sss{31}  & \sss{11}          & \sss{53}  & \sss{3}           & \sss{23}             & \sss{7}                   &  \sss{5}            &\sss{3,11}\\
\hline
\sss{d_p}      & \sss{13} & \sss{13} & \sss{20}           & \sss{9}   & \sss{14}          & \sss{10}  & \sss{12}          & \sss{14}  & \sss{21}          & \sss{12}  & \sss{50}          & \sss{16}             & \sss{29}                  &  \sss{36}           & \sss{55,43}\\
\hline
\end{tabular}
\end{center}
\end{table}
\medskip

From the above definitions it follows that for a fixed pair $(e,a)$ satisfying 
conditions~(C1) and~(C2) and such that
$d(e,p_a)\le ea$, inequality~\eqref{080211a} holds for all $d\ge ae$ and $p\in P(a,e)$. This 
enables us to find all the pairs for which inequality~\eqref{080211a} does not hold for at least one $p\in P(a,e)$ (the corresponding values of $a$ in Table~\ref{tbl1} are written in bold script):
\nmrt
\item[$\bullet$] $(e,a)=(13,3)$ or $(11,1)$,
\item[$\bullet$] $e=9$ and $a\in\{1,2,3\}$,
\item[$\bullet$] $e=5$ and $a\in\{1,\ldots,8\}$.
\enmrt
For each of these pairs we have to check inequality~\eqref{060211b} for all positive
integers $d\le d(e,3)$ which is a multiple of $ae$. The {\it available} triples $(e,a,d)$,
i.e. those that are obtained in this way, are listed in the first three rows of Table~\ref{tbl2}
below.\footnote{We did not cited in the table some available triples, like $(e,a,d)$ with 
$d\ge 22$ and $(e,a)=(11,1)$, because if the inequality \eqref{150211g} holds for some $d$, 
then it holds also for largest~$d$'s.} In the fourth 
and the fifth rows of this table we give respectively the values $p=p(a,e)$ and $q=q(e,a,d)$ 
where the latter number is equal to the minimal prime in $P(a,e)$ for which
\qtnl{150211g}
q^{d-\lfloor 4d/9\rfloor}>4\cdot ((q^a-1)/t_a\cdot e^2\cdot s_e\cdot a_0)^2;
\eqtn
Here the integer $t_a$ is defined as follows: if $a$ is odd, then $t_a=2$, otherwise 
$t_a=2^{t+2}$ where $t$ is the maximal positive integer such that~$2^t$ divides~$a$.
Then obviously $t_a$ divides $p^a-1$ for any odd prime~$p$, and hence
$u_{a,p}\le (p^a-1)/t_a$. Thus the required inequality~\eqref{060211b} follows
from~\eqref{150211g}. In the computation of $q$ we used the values of $s_e$ 
(and in cases $(e,a,d)=(5,3,15)$ and $(5,5,25)$ also the equality $|K:F|=3$) given in 
the second statement of Lemma~\ref{150211j}.\medskip
\begin{table}[bht]
\caption{}\label{tbl2}
\begin{center}
\begin{tabular}[tc]{|c|c|c|c| c|c|c|c| c|c|}
\hline
e & 13 & 11 &  9        & 9  & 9  & 9  & 5         & 5  & 5  \\
\hline
a & 3  & 1  &  1        & 1  & 3  & 2  & 4         & 4  & 8  \\
\hline
d & 39 & 11 &  \bf 9    & 18 & 27 & 18 & \bf 20    & 40 & 40 \\
\hline
p & 3  & 23 &  7        & 7  & 7  & 5  & 3         & 3  & 3  \\
\hline
q & 3  & 23 &  61       & 5  & 5  & 5  & 7         & 3  & 3  \\
\hline
\hline
e  &   5     &  5 & 5         & 5  & 5         & 5  & 5         & 5  & 5  \\
\hline
a  &   1     &  1 & 2         & 2  & 3         & 3  & 5         & 6  & 7  \\
\hline
d  &   \bf 5 & 10 & \bf 10    & 20 & \bf 15    & 30 &  25       & 30 & 35 \\
\hline
p  &   11    & 11 & 11        & 11 & 11        & 11 & 11        & 11 & 11 \\
\hline
q  &   5641  & 11 & 19        & 3  & 31        & 3  & 11        & 7  & 11 \\
\hline
\end{tabular}
\end{center}
\end{table}

It follows from the definition of $q$ that if $(e,a,d)$ is one of the available triples 
and $q\le p$, then inequality~\eqref{060211b} holds for all $p\in P(a,e)$. The remaining 
$5$ triples are the following: $(9,1,9)$, $(5,4,20)$, $(5,1,5)$, $(5,2,10)$ and $(5,3,15)$ 
(the corresponding values of $d$ in Table~\ref{tbl2} are written in bold 
script). For any of them the inequality~\eqref{060211b} does not hold only for 
those quadruples $(e,a,d,p)$ in which 
$$
p\in P(a,e)\cap\{1,\ldots,q-1\}.
$$ 
A straightforward 
check shows that these quadruples are exactly those listed in the lemma statement.\bull

\lmml{080411a}
The implication \eqref{220411a} holds unless
$(e,a,d)=(9,1,9)$ and $p\in\{7,19\}$, or $(e,a,d)=(5,1,5)$ and
$p\in\{11,31,41,61,71,101,151,181,271\}$.
\elmm
\proof  Given a prime $q$ denote by $k_q$ the number of all non-identity elements $g\in K$
the order of which is a power of~$q$; the maximum of $\fix(g)$ over all these elements~$g$ 
is denoted by~$f_q$. Clearly,
this maximum is achieved on the elements of order~$q$. The number $f_{q'}$ is defined in a
similar way: the maximum is taken over all non-identity elements $g\in K$ the order of
which is not a power of~$q$. Then it is easily seen that $\Fix(K)\le k_qf_q+(k-k_q)f_{q'}$.
So by Corollary~\ref{120211b} it suffices to prove that the inequality
\qtnl{090411a}
p^d>4(k-1)(k_q f_q+(k-k_q)f_{q'}).
\eqtn
holds for an appropriate prime divisor~$q$ of~$k=|K|$. By Lemma~\ref{050311a} it suffices to 
check this inequality only for those groups $K$ the parameters $(e,a,d,p)$ of which are listed 
in the statement of this lemma.\medskip

Let $(e,a,d)=(9,1,9)$ and $p\in\{13,31,37,43\}$. Then from Theorem~\ref{020311a} it follows
that $K=A$, $|F:U|=3^4$ and $U$ is a central subgroup of~$K$. Besides, by Lemma~\ref{150211j} 
we also have $|A:F|=s_e=5$. Thus
$$
k=|F|\cdot 5\quad\text{divides}\quad p_0:=\frac{p-1}{2}\cdot 3^4\cdot 5.
$$
It follows that the order of a Sylow $5$-subgroup of~$K$ is $5$ or $25$ depending on
whether $p\in\{13,37,43\}$ or $p=31$; in the former case $k_5=4\cdot 3^4$, whereas in
the latter one $k_5=20\cdot 3^4$. Moreover, in any case one can easily deduce from
Lemma~\ref{200411z} that $f_5\le p^{\lfloor d/5\rfloor}=p$ and $f_{5'}\le p^{d/3}=p^3$. 
Thus
$$
k_5 f_5+(k-k_5)f_{5'}\le k_5p+(p_0-k_5)p^3.
$$
A straightforward computation shows that the right-hand side of this inequality is less
than $p^9/4(k-1)$ for $p\in\{13,31,37,43\}$. This proves required inequality~\eqref{090411a} 
in our case.\medskip

A similar argument works when the quadruple $(e,a,d,p)$ is equal to $(5,4,20,3)$,
$(5,3,15,11)$ or $(5,2,10,11)$. In all these cases $|K:F|=3$ by Lemma~\ref{150211j}.
From now on we always assume that the group $U$ is the maximal odd subgroup of the 
multiplicative subgroup of $\spa(U)=\GF(p^a)$ (the base of a scheme is not decreased 
under taking a fusion). Then from Theorem~\ref{020311a} and Lemma~\ref{200411z} it follows that $k=75\cdot u$ with $u=|U|$, $f_3\le p^{\lfloor 4d/9\rfloor}$ and
\nmrt
\item[$\bullet$] if $(e,a,d,p)=(5,4,20,3\phantom{0})$, then $u=5$, $k_3=5^2\cdot 2$
and $f_{3'}\le p^4$,
\item[$\bullet$] if $(e,a,d,p)=(5,3,15,11)$, then $u=35\cdot 19$, $k_3\le 7\cdot 19\cdot 
5^2\cdot 2$ and $f_{3'}\le p^3$,
\item[$\bullet$] if $(e,a,d,p)=(5,2,10,11)$, then $u=15$, $k_3=5^2\cdot 6$ and $f_{3'}\le p^2$.
\enmrt
In all these cases inequality~\eqref{090411a} with $q=3$ follows by a straightforward
computation (for $d=15$ we have even more strong inequality in which the summand $(k-k_q)f_{q'}$
is replaced by $kf_{q'}$).\medskip

Let $(e,a,d)=(5,1,5)$ and $p$ belongs to the set $\P_0$ of all primes $q\le 5591$ such
that $q=1\,(\hspace{-6pt}\mod 5)$. In this case as before we have: $K=A$, $|A:F|=3$,
$|F:U|=5^2$ and $U$ is a central subgroup of $K$ (Theorem~\ref{020311a} and Lemma~\ref{150211j}).
Thus
\qtnl{230411ax}
k=|F|\cdot 3\quad\text{divides}\quad \frac{p-1}{2}\cdot 5^2\cdot 3.
\eqtn
A straightforward computation for all $p\in\P_0$ shows that the order of a Sylow $3$-subgroup
of the group $K$ equals to $3^{t_p}$ where $1\le t_p\le 6$. By Lemma~\ref{200411z} we have
\qtnl{240411y}
k_3=25(3^{t_p}-3^{t_p-1}),\quad f_3\le p^{\lfloor 4d/9\rfloor}=p^2,\quad
f_{3'}\le p^{\lfloor d/3\rfloor}=p
\eqtn
Next, given a positive integer $t$ denote by $p_0=p_0(t)$ the minimal prime in~$\P_0$ for which
$3^{t-1}$ divides~$(p_0-1)/2$, and by $p_1=p_1(t)$ the maximal real root of the polynomial
$$
g_t(x)=x^5-4\cdot(x'-1)\cdot( t'x^2+(x'-t')x).
$$
where $x'=75(x-1)/2$ and $t'=25(3^t-3^{t-1})$. When $t=t_p$ and $x=p$, we obtain
from~\eqref{230411ax} and~\eqref{240411y} that $x'\ge k$ and $t'=k_3$. Therefore 
$$
p^5\ge p^5-g_{t_p}(p)\ge 4\cdot (k-1)\cdot(k_3f_3+(k-k_3)f_{3'}).
$$
On the other hand, it is easily seen that $g_t(p)>0$ for all $t\ge 1$ and all $p>p_1(t)$.
Thus inequality~\eqref{090411a} does not hold only if $t\in\{1,\ldots,6\}$
and $p\in\P_0$ is such that $p_0(t)\le p\le p_1(t)\}$. In the Table~\ref{tbl3}
we present computed values of the functions $p_0(t)$ and $p_1(t)$. It follows
that the required inequality does not hold only if $(p,t)$ is one of the following
pairs: $(271,4)$, $(181,3)$, $\{31,61,151\}\times\{2\}$ and $\{11,41,71,101\}\times\{1\}$.
Thus the proof in this case is completely done.\bull
\begin{table}[bht]
\caption{}\label{tbl3}
\begin{center}
\begin{tabular}{|c|c|c|c|c|c|c|}
\hline
$t$                   & 1   & 2   &  3  & 4   & 5   & 6\\
\hline
$p_0$                 & 11  & 31  & 181 & 271 & 811 & 4861\\
\hline
$\lfloor p_1\rfloor$  & 113 & 166 & 269 & 455 & 782 & 1351\\
\hline
\end{tabular}
\end{center}
\end{table}

\lmml{140411a}
The implication \eqref{220411a} holds for all $e>1$.
\elmm
\proof  We recall that the set $\Omega$ is identified with a $d$-dimensional linear space 
over a field $\GF(p)$ and the group $G$ contains the translation group of~$\Omega$. Denote 
by $\alpha$ the zero vector of~$\Omega$. Then given $r\in S$ and $\beta\in\alpha r$ 
the intersection number $c_{ts}^r$ is equal to the number of all $\gamma\in\alpha t$
for which $r(\gamma,\beta)=s$. Besides, it is easily seen that $r(\gamma,\beta)=r(\gamma',\beta)$ if and only if $\gamma-\beta\in(\gamma'-\beta)^K$ for all $\gamma'$. Thus
\qtnl{030811a}
c_{ts}^r=|\Delta_\beta(t)|
\eqtn
where $\Delta_\beta(t)$ the set of all sets $(\gamma-\beta)^K\cap(\alpha t-\beta)$
with $\gamma\in\alpha t$ and $\alpha t-\beta$ is the set of all vectors $\gamma'-\beta$ 
with $\gamma'\in\alpha t$. Then using~\eqref{030811a} for $t=r$ and $t=r^*$ we can compute 
the numbers $\rsc{r}{2}{r}$ and $|\rsc{r}{2}{r^*}|$ defined before Lemma~\ref{160411a} as follows:
$$
|\rsc{r}{2}{r}|=|\{\Delta\in\Delta_\beta(r):\ |\Delta|\le 2\}|=:a_r
$$
and
$$
|\rsc{r}{2}{r^*}|=|\{\Delta\in\Delta_\beta(r^*):\ |\Delta|\le 2\}|=:b_r
$$
(in the second case we used equalities $c_{rs}^{r^*}=c_{s^*r^*}^r=c_{r^*s^*}^r$ that
follow from~\eqref{150410a}). Our goal is to find a relation $r\in S$ such that
\qtnl{190411a}
a_r+b_r>2n_r/3.
\eqtn
Then Lemma~\ref{160411a} implies that $b(\X)\le 2$, and we are done. To find such $r$ 
we can assume that $(e,a,d,p)$ is one of the quadruples listed in the statement
of Lemma~\ref{080411a}.\medskip

Suppose first that $(d,p)=(9,7)$. Denote by $E$ an extraspecial group of order~$3^5$ 
and exponent~$3$. Then $K=A$ is isomorphic to a semidirect product $K_0=E.5$ 
in which the group of order $5$ acts irreducibly on $E/Z(E)$. By means of GAP we found
that $K_0$ is uniquely determined up to isomorphism (there is a unique non-nilpotent
group of order $3^5\cdot 5$ with a nonabelian Sylow $3$-subgroup). Moreover, up to equivalence 
there are exactly two classes of irreducible $d$-dimensional $K_0$-modules over $\GF(p)$. For 
both of them we constructed in GAP generators for the corresponding primitive subgroup 
of $\GL(d,p)$ isomorphic to~$K$ (see Section~\ref{240411q}). Then we fixed a standard linear 
base $\{e_1,\ldots,e_d\}$ in $\GF(p)^d$ and took 
$$
\beta=e_1+e_2+e_5\quad\text{and}\quad r=r(\alpha,\beta).
$$
A straightforward computation shows that in both cases $n_r=|\beta^K|=|K|=1215$, $a_r=1035$ and $b_r\ge 754$. Therefore inequality~\eqref{190411a} do hold and we are done.\medskip

The computation in each of the remaining case is essentially the same as in the above
case $(d,p)=(9,7)$. The minor differences are the following. In the case $(d,p)=(9,19)$
we have $K=\lg K_0,\xi_p I_d\rg$ where $I_d$ is the identity matrix in $\GL(d,p)$,
and $\xi_p\in\GF(p)$ is a generator of the maximal multiplicative $2'$-subgroup in~$\GF(p)$
(in our case, the subgroup of order~$9$). In the case $(5,p)$ the group $E$ is an
extraspecial group of order $5^3$ and exponent~$5$, $K_0$ is a semidirect product $E.3$ 
in which the group of order $3$ acts irreducibly on $E/Z(E)$, and $K=\lg K_0,\xi_p I_d\rg$.
The computation results cited in the Table~\ref{tbl4} below
\begin{table}[bht]
\caption{}\label{tbl4}
\begin{center}
\begin{tabular}{|c|c|c|c|c|c|c|c|c|c|c|c|}
\hline
   & \multicolumn{2}{c|}{$(9,9,1)$} & \multicolumn{9}{c|}{$(5,5,1)$}\\
\hline
$p$         & 7       & 19   & 11  & 31    & 41  & 61   & 71   & 101  & 151  & 181  & 271\\
\hline
$\beta$     & \multicolumn{2}{c|}{$e_1+e_2+e_5$} &  \multicolumn{9}{c|}{$e_1+e_2$} \\
\hline
$n_r$       & 1215    & 3645 & 375 & 1125  & 375 & 1125 & 2625 & 1875 & 5625 & 3375 & 10125 \\
\hline
$a_r$       & 1035    & 3483 & 199 &  987  & 361 & 1061 & 2413 & 1755 & 5221 & 3199 & 9421 \\
\hline
$b_r$       &  754    & 1697 &  99 &  469  & 160 &  526 & 1181 &  853 & 2585 & 1563 & 4685 \\
\hline
$N$         & 2       &   2  &  4  & 12    &  4  &  12  &  4   &  4   & 12   & 12   & 12\\
\hline
\end{tabular}
\end{center}
\end{table}
show that inequality~\eqref{190411a} holds in all the cases (in the last row in the table  
we gives the number of the classes of irreducible $d$-dimensional $K_0$-modules over $\GF(p)$; 
the  values $a_r$ and $b_r$ correspond to the $K_0$-module with minimal sum $a_r+b_r$).\bull

{\bf Proof of Theorem~\ref{291210d}.}
Let $\X=\inv(G)$ where $G\le\AGL(d,p)$ is an affine group with primitive zero 
stabilizer~$K\le\GL(d,p)$ of odd order, $p$ is an odd prime. Then from Lemmas~\ref{140411a}
and~\ref{150211j} it follows that $K$ is contained in the unique Hall $2'$-subgroup $K^*$
of the group $\GaL(1,p^d)$. It is easily seen that $\orb(K^*,\Omega)=\orb(K',\Omega)$
where $K'$ is the maximal odd order multiplicative group of the field
$\FF=\GF(p^d)$. Thus
\qtnl{300711b}
\X=\inv(G)\ge\inv(G^*)=\inv(G')=:\X'
\eqtn
where $G^*=AK^*$ and $G'=AK'$ with $A$ being the translation group of the linear
space~$\Omega$. However, $\X'$ is a cyclotomic scheme over the field~$\FF$. Therefore
$b(\X')\le 3$ by Theorem~\ref{200411b}. Thus by~\eqref{300711b}
we conclude that $b(\X)\le b(\X')\le 3$ which completes the proof.\bul

\section{The proof of Theorems~\ref{290311a}, \ref{200411c} and~\ref{140211x}}\label{050711u}

\sbsnt{Proof of Theorem~\ref{200411c}.}
We argue by induction on the degree of a schurian antisymmetric coherent 
configuration~$\X=\inv(G)$.
By the first inequality in~\eqref{200411e}, we can assume that it is homogeneous.
Suppose that $\X$ is imprimitive. Then there is a nontrivial equivalence 
relation~$e\in S^\cup$. The schemes $\X_1=\X_\Gamma$ where $\Gamma\in\Omega/e$,
and $\X_2=\X_{\Omega/e}$ are obviously antisymmetric and schurian. Moreover,
$\X$ is isomorphic to a fission of the scheme $\X_1\wr\X_2$. By Corollary~\ref{190311a}
this implies that
$$
gb(\X)\le gb(\X_1\wr\X_2)\le\max\{gb(\X_1),gb(\X_2)\}
$$
and we are done by induction.\medskip

Suppose that the scheme $\X$ is primitive. Then by Theorem~\ref{140311c}
it is either linearly imprimitive or linearly primitive. In the former case
$\X$ has a nontrivial fusion isomorphic to $\Y\pwr L$ where $\Y$ is a linearly primitive
antisymmetric scheme and $L$ is a transitive group of odd order (Theorem~\ref{160311a}). Thus
by induction and Theorem~\ref{140311d} we have
$$
gb(\X)\le gb(\Y\pwr L)\le\max\{gb(\Y),gb(\inv(L))\}\le 1.
$$
To complete the proof suppose that $\X$ is linearly primitive. If it is
cyclotomic, then we are done by Theorem~\ref{200411b}. Otherwise by
Theorem~\ref{291210d} it has a base $B=\{\alpha,\beta\}$ where $\alpha$ and $\beta$
are two (possibly equal) points in~$\Omega$. In this case 
$$
\fis(\X,\{B\})=\X_{\alpha,\beta}
$$
because the scheme $\X$ is antisymmetric. Thus $gb(\X)\le1$ and we are done.\bull

\sbsnt{Proof of Theorem~\ref{290311a}.}
By Theorem~\ref{140311c} the scheme $\X=\inv(G)$ is either linearly imprimitive or 
linearly primitive. In the latter case we are done by Theorem~\ref{291210d}.
In the former case $\X$ has a nontrivial fusion isomorphic to $\Y\pwr L$ where $\Y$ is
a linearly primitive antisymmetric scheme and $L$ is a transitive group of odd order
(Theorem~\ref{160311a}). In particular, the scheme $\Y$ is not complete. Besides,  
$b(\Y)\le 3$ by Theorem~\ref{291210d} and $b:=gb(\inv(L))=1$ by Theorem~\ref{200411c}. 
This implies by Theorem~\ref{140311d} that
$$
b(\X)\le b(\Y\pwr L)\le  b(\Y)+\max\{0,1-\lceil (b(\Y)-1)/2\rceil\}.
$$
When $b(\Y)=1,2,3$, the right-hand side of the above inequality is equal respectively to $2,2,3$. In any case $b(\X)\le 3$, and we are done.\bull

\sbsnt{Algorithm.}
We will deduce Theorem~\ref{140211x} from Theorem~\ref{270411a} proved below.
The algorithm constructed in the proof of the latter theorem is, in a sense, a combinatorial
version of the Babai-Luks algorithm from~\cite{BL83}. The following statement to be used
in the proof of Theorem~\ref{270411a} is a special case of Corollary~3.6 of that paper.
In what follows we always assume that any permutation group on the input or output of
an algorithm is given by a generator set of polynomial size in the degree of the group.

\thrml{blca}
Let $G\le\sym(\Omega)$ be a solvable group of degree~$n$. Then given a coherent configuration
$\X$ on $\Omega$, the group $\aut(\X)\cap G$ can be found in time $n^{O(1)}$.\bul
\ethrm

To apply Theorem~\ref{blca} we have to be able to construct the group~$G$. This will be done by means
of Corollary~\ref{290411a} and the following statement proved in \cite[Lemma~3.4]{EP03be}.
Below for permutation groups $G_1,\ldots,G_s$, $s\ge 1$, we define the group
$\Wr(G_1,\ldots,G_s)$ to be the iterated wreath product 
$(\cdots(G_1\wr G_2)\wr\cdots)\wr G_s$ in imprimitive action. 

\lmml{theorem4}
Let $\X$ be a scheme and $1_\Omega=e_0\subset e_1\subset\ldots\subset e_s=\Omega^2$ a series
of equivalence relations of~$\X$. Suppose that for $i=1,\ldots s$ a permutation group $G_i$ on a
set~$\Delta_i$ and a family of bijections $f_\Gamma:\Gamma\to\Delta_i$ where
$\Gamma\in\Omega_i$ with $\Omega_i=\{\Gamma'/e_{i-1}:\ \Gamma'\in\Omega/e_i\}$, form a 
majorant of $\aut(\X)$ with respect to the pair~$(e_{i-1},e_i)$. Then the mapping 
$$
f:\Omega\to\prod_{i=1}^s\Delta_i,\quad \alpha\mapsto(\ldots,f_i(\Gamma_{i-1}),\ldots)
$$
is a bijection and $\aut(\X)^f\le\Wr(G_1,\ldots,G_s)$ where $f_i=f_{\Gamma_i}$ and 
$\Gamma_{i-1}$ and $\Gamma_i$ are respectively the classes of~$e_{i-1}$ and~$e_i$ containing $\alpha$.\bull
\elmm

{\bf Proof of Theorem~\ref{270411a}.}
To describe the algorithm we need the auxiliary procedure $\Test(\X,G)$ that 
given a coherent configuration $\X=(\Omega,S)$ and a group $G\le\sym(\Omega)$ output 
$G$ or empty set depending whether or not $\X=\inv(G)$.
Since the latter equality exactly means that $S=\orb(G,\Omega^2)$, the procedure
can be implemented in polynomial time in $|\Omega|$ by means of a standard algorithm
finding the orbits of a permutation group (see e.g.~\cite{S02}).\medskip\medskip

\centerline{\bf Schurity Recognition Algorithm}\medskip

\noindent {\bf Input:} an antisymmetric coherent configuration $\X$.

\noindent {\bf Output:} the group $\aut(\X)$, or $\X$ is not schurian.\medskip

\noindent{\bf Step 1.}
If $\X$ is not homogeneous, then recursively apply the algorithm to
the coherent configuration $\X_1=\X_{\Delta_1}$ and $\X_2=\X_{\Delta_2}$ 
where $\Delta_1$ is a fiber of~$\X$ and $\Delta_2$ is its complement.
If either $\X_1$ or $\X_2$ is not schurian, then so is~$\X$, else output
$\Test(\X,H)$ where $H\le\sym(\Omega)$ is the group found 
by the algorithm of Theorem~\ref{blca} for $G=\aut(\X_1)\times\aut(\X_2)$.\medskip

\noindent{\bf Step 2.}
Find a maximal series of equivalence relations as in the
hypothesis of Lemma~\ref{theorem4}. If there exist $i\in\{1,\ldots,s\}$ and
$\Gamma,\Gamma'\in\Omega_i$ such that
\qtnl{040511a}
b(\X_\Gamma)>3\quad\text{or}\quad\iso(\X_{\Gamma^{}},\X_{\Gamma'},\varphi_{\SG,\SG'})=\emptyset
\eqtn
where $\varphi_{\SG,\SG'}$ is 
the algebraic isomorphism~\eqref{140711a} for $\X=\X_{\Omega/e_{i-1}}$, then
$\X$ is not schurian.\medskip

\noindent{\bf Step 3.}
By the algorithm of Corollary~\ref{290411a} find a majorant of $\aut(\X)$ with respect 
to $(e_{i-1},e_i)$, say $G_i\le\sym(\Delta_i)$ and $\{f_\Gamma\}_{\Gamma\in \Omega_i}$
for $i=1,\ldots s$.\medskip

\noindent{\bf Step 4.}
Output $\Test(\X,H)$ where $H=\aut(\X)\cap G^{f^{-1}}$ is the group
found by the algorithm of Theorem~\ref{blca} with $G=\Wr(G_1,\ldots,G_s)$ and $f$ as in
Lemma~\ref{theorem4}.\bul\medskip

To prove the correctness of the algorithm suppose first that the coherent configuration~$\X$
is not homogeneous. Then it is schurian only if so are the coherent configurations~$\X_1$
and~$\X_2$ found at Step~1, and, moreover,
$$
\aut(\X)\le\aut(\X_1)\times\aut(\X_2)
$$
where the group in the right-hand side has odd order. Thus the correctness in this case
follows from Theorem~\ref{blca}. Let $\X$ be homogeneous. Then it is schurian only if
for all $\Gamma $ and $\Gamma'$ defined at Step~2 we have
$$
\X_\Gamma=\inv(\aut(\X)^\Gamma)\qaq
\aut(\X)_{\Gamma\to\Gamma'}\subset\iso(\X_{\Gamma^{}},\X_{\Gamma'},\varphi_{\SG,\SG'})
$$
where $\aut(\X)_{\Gamma\to\Gamma'}$
is the set of all bijections from $\Gamma$ onto $\Gamma'$ induced by the automorphisms 
of~$\X$. On the other hand, the maximality condition in choosing~$e_i$'s implies that 
under the schurity assumption each scheme $\X_\Gamma$ is also primitive. Therefore 
by Theorem~\ref{290311a} in our case $b(\X_\Gamma)\le 3$ for all~$\Gamma$. Thus the
relations~\eqref{040511a} imply that~$\X$ is not schurian, and the output of Step~3 is
correct. Now, the correctness of the output at Step~4 follows from Lemma~\ref{theorem4}.
Finally a polynomial bound for the running time of the algorithm follows from
Theorem~\ref{blca} and Corollary~\ref{290411a}.\bull

\sbsnt{Proof of Theorem~\ref{140211x}.}
Given a colored tournament $T$ denote by $\X=\X(T)$ the coherent configuration
$\X(T)=\fis(\cS)$ where $\cS$ is the set of color classes of the arc set 
of~$T$. Then $T$ is schurian if and only if the coherent configuration~$\X$ is schurian. 
Since the latter can be constructed in time $n^{O(1)}$ where $n$ is the number of
vertices of~$T$ (see Subsection~\ref{310711a}), statements~(1) and~(2) immediately follow
from Theorem~\ref{270411a}.\medskip 

To prove statement~(3) let $T_i$ be a colored schurian tournament, $\cS_i$ the set of 
color classes of the arc set of~$T_i$ and $\X_i=\X(T_i)$, $i=1,2$. Then by
Theorem~\ref{thbw} without loss of generality we can assume that there exists an algebraic
isomorphism $\varphi:\X_1\to\X_2$ such that
\qtnl{040511c}
\cS_1^\varphi=\cS_2
\eqtn
(for otherwise $\iso(T_1,T_2)=\emptyset$). In this case $\iso(T_1,T_2)=\iso(\X_1,\X_2,\varphi)$.
To construct the latter set take a copy $\X_3$ of the coherent configuration $\X_2$.
Set
$$
\W=\{\X_i\}_{i=1}^3\quad\text{and}\quad\Psi=\{\psi_{i,j}\}_{i,j=1}^3
$$
where $\psi_{i,j}:\X_i\to\X_j$ is an algebraic isomorphism defined as follows (below $S_i$ denotes the set of basis relations of~$\X_i$): 
\nmrt
\item[$\bullet$] $\psi_{i,j}=\id_{S_i}$, if $i=j$ or $\{i,j\}=\{2,3\}$, 
\item[$\bullet$] $\psi_{i,j}=\varphi$,\quad\ if $1=i\ne j$,
\item[$\bullet$] $\psi_{i,j}=\varphi^{-1}$, if $i\ne j=1$.
\enmrt
 According to
\cite[Definition~7.1]{EP99c} there exists the smallest coherent configuration~$\X$
on the disjoint union $\Omega=\Omega_1\cup\Omega_2\cup\Omega_3$ where $\Omega_i$ is the
point set of~$\X_i$, such that
$$
\X_{\Omega_i}=\X_i\quad\text{and}\quad\X_{\Omega/e}=\inv(G)
$$
where $i=1,2,3$, $e=\Omega_1^2\cup\Omega_2^2\cup\Omega_3^2$ and $G$ is the cyclic subgroup
of $\sym(3)$. It was also proved there (see \cite[Corollary~7.9]{EP99c}) that $\X$ is
schurian if and only if $\X_i$ is schurian for all~$i$ and the algebraic
isomorphism $\psi_{i,j}$ is induced by an isomorphism for all~$i,j$. In our case
the former condition is obviously satisfied, whereas the latter one is satisfied
if and only if the set $\iso(\X_1,\X_2,\varphi)$ is not empty.\medskip

To complete the proof we note that the coherent configuration~$\X$ is
antisymmetric. Therefore by Theorem~\ref{270411a} one can test whether or not~$\X$ is schurian and (if so) find the group $\aut(\X)$ in time $n^{O(1)}$. Now, if $\X$ is not schurian, then by the above 
$$
\iso(\X_1,\X_2,\varphi)=\emptyset.
$$ 
On the other hand, if $\X$ is schurian, then by means of standard permutation
group algorithms (see e.g.~\cite{S02}) one can efficiently find an element
$g\in\aut(\X)$ such that $\Omega_1^g=\Omega_2$, and the setwise stabilizer 
$H$ of the set $\Omega_1$ in the group $\aut(\X)$.
Since in this case obviously 
$$
\iso(\X_1,\X_2,\varphi)=\{h^{\Omega_1}g_{\Omega_1}:\ h\in H\}
$$
where $h^{\Omega_1}$ is the restriction of $h$ on $\Omega_1$, and 
$g_{\Omega_1}:\Omega_1\to\Omega_2$ is the bijection induced by~$g$, we are done.\bull


\begin{thebibliography}{99}

\bibitem{ADM}
V.~Arvind, B.~Das, P.~Mukhopadhyay, {\em Isomorphism and canonization of tournaments and
hypertournaments}, Journal of Computer and System Sciences, {\bf 76} (2010), 509–-523.

\bibitem{B81}
L.~Babai, {\em On the order of uniprimitive permutation groups}, Annals of Math.,
{\bf 113} (1981), 553--568.

\bibitem{BL83}
L.~Babai, and E.~M.~Luks, {\em Canonical labeling of graphs}, Proc. 15th ACM
STOC, (1983), 171--183.

\bibitem{CB}
R.~F.~Bailey, P.~J.~Cameron, {\em Base size, metric dimension and other invariants
of groups and graphs}, Bull. LondonMath. Soc., {\bf 43} (2011) 209-–242.

\bibitem{EP99}
S.~Evdokimov, I.~Ponomarenko, {\em On primitive cellular algebras}, Zapiski Nauchnykh Seminarov
POMI, {\bf 256} (1999), 38--68. English translation in J. Math. Sci. (New York),
{\bf 107} (2001), no.5, 4172--4191.

\bibitem{EP99c}
S.~Evdokimov, I.~Ponomarenko, {\em On highly closed cellular
algebras and highly closed isomorphisms}, Electronic J. of Combinatorics,
{\bf 6} (1999), \#R18.

\bibitem{EP01}
S.~Evdokimov, I.~Ponomarenko, {\em Two-closure of odd permutation group in polynomial
time}, Discrete Mathematics, {\bf 235/1-3}, (2001), 221--232.

\bibitem{EP01ce}
S.~Evdokimov, I.~Ponomarenko, {\em Characterization of cyclotomic schemes and normal
Schur rings over a cyclic group}, Algebra and Analysis, {\bf 14} (2002), 2, 11--55.
English translation in  St. Petersburg Math. J.,  {\bf 14}  (2003),  no. 2, 189--221.

\bibitem{EP03be}
S.~Evdokimov, I.~Ponomarenko, {\em Recognizing and isomorphism testing circulant graphs
in polynomial time}, Algebra and Analysis, {\bf 15} (2003), 6, 1--34. English  translation
in  St. Petersburg Math. J., {\bf  15} (2004), no. 6, 813--835.

\bibitem{EP09}
S.~Evdokimov, I.~Ponomarenko, {\em Permutation group approach to association schemes},
European J. Combin., {\bf 30} (2009), 6, 1456--1476.

\bibitem{EPT}
S.~Evdokimov, I.~Ponomarenko, G.~Tinhofer, {\em  Forestal
algebras and algebraic forests (on a new class of weakly compact graphs)},
Discrete Mathematics, {\bf 225} (2000), 149--172.

\bibitem{Es91a}
A.~Espuelas, {\em Regular orbits on symplectic modules}, J. Algebra, {\bf 138} (1991), 1--12.

\bibitem{Es91b}
A.~Espuelas, {\it Large character degrees of groups of odd order}, Illinois J. Math.,
{\bf 35} (1991), 499--505.

\bibitem{G83}
D.~Gluck, {\it Trivial set-stabilizers in finite permutation groups}, Can. J. Math.,
{\bf XXXV} (1983), 59--67.

\bibitem{JKM}
G.~A.~Jones, M.~Klin, Y.~Moshe,
{\em Primitivity of Permutation Groups, Coherent Algebras and Matrices},
Journal of Combinatorial Theory, {\bf A98} (2002), 210–-217.

\bibitem{GAP}
The GAP Group, GAP-4—Groups, Algorithms and Programming, Version 4.4.5, 2005,
{\tt http://www.gap-system.org}.

\bibitem{MW}
O.~Manz, T.~Wolf, {\it Representations of solvable groups}, London Math. Soc. Lecture
Note Ser. vol. 185, Cambridge Univ. Press, Cambridge (1993).

\bibitem{MP09}
 M.~Muzychuk, I.~Ponomarenko, {\em On Pseudocyclic Association Schemes},
{\tt arXiv:math.CO/0910.0682} (2009), 1--23.

\bibitem{S96}
A.~Seress, {\em The minimal base size of primitive solvable permutation groups}, J. London
Math. Soc., {\bf 53} (1996) 243–-255.

\bibitem{S02}
A.~Seress, {\em Permutation Group Algorithms}, Cambridge Univ. Press, 2002.


\bibitem{Su76}
D.~A.~Suprunenko, {\em Matrix Groups}, American Mathematical Society, Providence, RI, 1976.

\bibitem{YY10}
Yong Yang, {\em Regular orbits of finite primitive solvable groups},
Journal of Algebra, {\bf 323} (2010), 2735-–2755.

\bibitem{W76}
B.~Weisfeiler (editor), {\em On construction and identification of graphs}, Springer Lecture
Notes, 558, 1976.

\bibitem{Zi1}
P.-H.~Zieschang, {\em Theory of Association Schemes}, Springer, Berlin \& Heidelberg,
(2005).

\end{thebibliography}
\end{document}